\documentclass[10pt,a4paper]{article}
\usepackage{latexsym,amsfonts,amssymb,amsthm, mathrsfs}
\usepackage{colordvi,multicol}
 \usepackage[leqno]{amsmath}
\usepackage{cases}
\usepackage{color}
\newcommand{\ooo}{{\overline{\bf E}}}
\newcommand{\uuu}{\underline {\bf E}}

\newtheorem{thm}{Theorem}%%[section]
\newtheorem{cor}{Corollary}
\newtheorem{lem}{Lemma}

\newtheorem{rem}{Remark}

\makeatletter
\def\@makefnmark{}
\makeatother

\begin{document}
\title{\bf   Prokhorov distance
\\ with rates of convergence
\\ under   sublinear expectations${}^{1)}$ }
\author{ \small QIANQIAN ZHOU$^*$,
\small ALEXANDER SAKHANENKO$^{**}$    \footnotesize{ AND}
\small JUNYI  GUO$^{***}$
 }

\date{}
\maketitle

\vskip 0.5cm \noindent{\bf Abstract}\quad
Prokhorov distances under sublinear expectations are presented in CLT and functional CLT,  and the convergence rates for them are obtained by Lindeberg method. In particular, the obtained estimate in functional CLT yields known Borovkov's estimate in classical  functional CLT with explicit constant.
\smallskip


\smallskip

\footnote{${}^*$ School of Mathematical Sciences,
 Nankai University,  Tianjin 300071,  P.R. China.\\ Email: qianqzhou@yeah.net.

 ${}^{**}$ Sobolev Institute of Mathematics,
4 Acad. Koptyug avenue,  Novosibirsk 630090, Russia.

\hspace{4mm}School of Mathematical Sciences,
 Nankai University,  Tianjin 300071,  P.R. China.
\\ Email: aisakh@mail.ru.

 ${}^{***}$ School of Mathematical Sciences,
 Nankai University,  Tianjin 300071,  P.R. China.\\  Email: jyguo@nankai.edu.cn.

${}^{1)}$This work was supported by the National Natural Science
Foundation of China (Grant No. 11571189) and the Fundamental Research Fund for
the Central Universities.

The research of A.Sakhanenko  was also supported by the program of fundamental scientific researches of the SB RAS № I.1.3., project № 0314-2019-0008

% Postal address: School of Mathematical Sciences,
% Nankai University, Tianjin 300071, P.R. China.\\
% \indent$^{\dagger}$Email: qianqzhou@yeah.net,
% \indent$^{\ddagger}$Email:��aisakh@mail.ru,
%\indent$^{*}$Email: jyguo@nankai.edu.cn.
}

\section{Introduction }

{\bf 1.1.} Convergence of a sequence of  random variables is an important concept in probability theory and its applications.
Convergence in distribution is the weakest form of such convergences. However, convergence in distribution is very frequently used in practice, one of the most important  application is  the central limit theorem (CLT).
It is well known that if we have classical i.i.d.  random variables $X_{1},X_2, \ldots$ with
\begin{gather}                                                                    \label{a1}
{\bf E}[X_1]=\mu, \quad 0<\sigma^{2}={\bf Var} X_{1}={\bf E}[(X_{1}-\mu)^{2}]<\infty,
\end{gather}
then (all limits in the paper are taken as $n\longrightarrow \infty$)
\begin{gather}                                                                    \label{a2}
W_n:=\frac{1}{\sqrt{n}}\sum_{i=1}^{n}\frac{X_i-\mu}{\sigma}\Rightarrow Z\sim N(0,1),
\end{gather}
where $Z\sim N(0,1)$ means that  $Z$ has the standard normal distribution.
In particular, convergence (\ref{a2})  is equivalent to the fact that for all bounded and  continuous function $\varphi(\cdot)$ we have
\begin{gather}                                                                    \label{a3}
{\bf E}[\varphi(W_n)]\longrightarrow {\bf E}[\varphi(Z)].
\end{gather}
And (\ref{a2})  %%the letter convergence
is also equivalent to the convergence
$$\pi (W_n, Z)\rightarrow0,$$
where $\pi (W_n, Z)$ is the Prokhorov distance between the distributions in $\mathbb{R}$ of random variables $W_n$ and $Z.$
We may define this Prokhorov distance in the following way:
\begin{eqnarray*}
\pi(W_n, Z)=\inf\big\{\varepsilon>0:\sup\big[{\bf P}(W_n \in A)-{\bf P}(Z \in A^{\varepsilon})\big] \leq \varepsilon\big\},
\end{eqnarray*}
where the supremum is taken over all measurable sets $A\subset \mathbb{R}$,  and $A^{\varepsilon}$ denotes the $\varepsilon$-neighborhood of $A,$ i.e., \
%\begin{gather*}         %%\label{a6}
$A^{\varepsilon}=\big\{x\in \mathbb{R}: \exists y\in A\quad\text{with} \quad |x-y|<\varepsilon\big\}.$

{\bf 1.2.}
The convergence in distribution of a sequence of random processes
in the functional central limit theorem has similar properties.
Consider again classical i.i.d. random variables  %%satisfying
from (\ref{a1}) and introduce into consideration the classical random broken line $S(\cdot)$  with the values
\begin{eqnarray}\label{1.8}
S(0)=0, \ \ S(k)=\sum_{i=1}^{k}\xi_i,\quad \xi_i:=\frac{X_{i}-\mu}{\sigma},\quad  i,k=1, 2, \ldots .\nonumber
\end{eqnarray}
Thus,
 $S(t)=S([t])+(t-[t])\xi_{[t]+1}$ %%\frac{X_{[t]+1}-\mu}{\sigma}$$
 is a linear function on each interval  of the form $[k, k+1],\  \  k=0,1,2,\ldots.$ Also introduce the sequence of rescaled broken lines
%%\begin{eqnarray}\label{1.9}
$S_{n}(t)={S(nt)}/{\sqrt{n}}$, %\quad
$t\in[0,1]$, %\quad
$n=1, 2, \ldots .$
%\nonumber%%\end{eqnarray}
Hence, $S_n(\cdot)$ is a broken line with the following values at key points:
\begin{gather}                                                 \label{a11}
S_n({k}/{n})=\frac{1}{\sqrt{n}}\sum_{i=1}^{k}\xi_i,\quad \xi_i=\frac{X_i-\mu}{\sigma}, \quad i,k=1,2, \ldots, n.
\end{gather}
We treat $S_n=S_n(\cdot)$ as a random process with values in the space  $\mathbb{C}[0,1]$ of all continuous functions, defined on the interval  $[0,1]$, with  the uniform norm~\mbox{$\|\cdot\|$}.
Then the sequence $S_n=S_n(\cdot)$ satisfies  functional  CLT,
introduced in \cite{D51} and~\cite{P56}, i.e., for any real bounded and continuous   functional $f(\cdot)$, defined on  $\mathbb{C}[0,1]$,  we have
\begin{gather*}         %                                        \label{a12}
{\bf E}[f(S_n)]\longrightarrow {\bf E}[f(B)],
\end{gather*}
where $B=B(\cdot)$ is a Wiener process.
This fact is also equivalent to the following convergence
\begin{gather*}      %                                           \label{a13}
\Pi(S_n, B)\longrightarrow 0,
\end{gather*}
where $\Pi(S_n, B)$ is the Prokhorov distance between the distributions in $\mathbb{C}[0,1]$ of the processes $S_n$ and $B$, i.e.,
\begin{eqnarray*}
\Pi(S_n, B)=\inf\big\{\varepsilon>0: %%\Pi(S_n, B, \varepsilon)
\sup\big[{\bf P}(S_n \in D)-{\bf P}(B \in D^{\varepsilon})\big]
\leq \varepsilon\big\},
\end{eqnarray*}
where the supremum is taken over all measurable sets $D\subset \mathbb{C}[0,1]$,  and $D^{\varepsilon}$ denotes the $\varepsilon$-neighborhood of $D,$ i.e.,
\begin{gather}\label{a16}
D^{\varepsilon}=\big\{x\in \mathbb{C}[0,1]: \exists y\in D\quad\text{with} \quad \|x-y\|<\varepsilon\big\},
\end{gather}
where $\|x-y\|=\max\{|x(t)-y(t)|:t\in[0,1]\}$ is the distance between $x,y\in \mathbb{C}[0,1]$. %% and $y\in \mathbb{C}[0,1].$

{\bf 1.3.}
Prokhorov distance was introduced in the seminal paper \cite{P56} of Prokhorov as an
useful instrument for investigation of convergence of probability distributions  in separable metric spaces. It
is also an  useful instrument for investigations of  rates of convergence in
functional CLT.
Indeed, if the functional $f(\cdot)$ is Lipschitz continuous, i.e.,
\begin{gather}\label{a17}
\forall x, y \in \mathbb{C}[0,1] \quad \big|f(x)-f(y)\big|\leq L\big\|x-y\big\|, \quad L<\infty,
\end{gather}
%%where $L$ is a Lipschitz constant,
 and if $f(B)$ has a bounded density, i.e.,
\begin{gather}\label{a18}
 \forall a\in \mathbb{R}\ \ \forall h\geq 0 \quad {\bf P}\big(a\leq f(B)\leq a+h\big)\leq K h,\quad K<\infty,
\end{gather}
%%where $K<\infty$ is a constant,
then (see \cite {B73} )
\begin{gather}                                                                                                                 \label{a19}
\sup_{a}\big|{\bf P} (f(S_n)\leq a)-{\bf P}(f(B)\leq a)\big|\leq (KL+1)\Pi(S_n, B).
\end{gather}

Inspired by inequality (\ref{a19}), a difficult task of obtaining estimates for the Prokhorov distance
$\Pi(S_n, B)$ attracted a lot of attention. The first such estimate, when ${\bf E}|\xi_1|^3<\infty,$
was obtained by  Prokhorov in  \cite{P56}. The first general estimate due to Borovkov
\cite{B73}. For i.i.d. random variables $\{\xi_i\}$ he proved in \cite{B73} that
\begin{gather}                                                                                                              \label{a41}
\Pi(S_n, B)\le C_1(p)\frac{
\big({\bf E}|\xi_1|^p\big)^{1/(p+1)}
}{n^{(p-2)/(2p+2)}},\quad p\in[2,3],
\end{gather}
where $C_1(p)$ depends only on $p$ and $\sup_{p\in[2,3]}C_1(p)<\infty$ is an  absolute  constant.
Later Sakhanenko in \cite{S89}  extended (\ref{a41}) for all $p\ge2$.
In \cite{S89} and \cite{S06} you may find  simple proofs of the fact that inequality (\ref{a41}) is unimprovable for all $p\ge2$, up to a constant $C_1(p)$.

{\bf 1.4.} In the above statements  the probability ${\bf P}$ and  expectation ${\bf E}$  are considered in  classical sense.
Recently, motivated by the risk measures, superhedge pricing and modeling uncertain in finance, Peng  in \cite{Pe05} (see, also \cite{Pe07b} and \cite{Pe10}) introduces  the notion of sublinear expectation. Thus,
 instead of ${\bf P}$ and ${\bf E}$ in the classical sense, we consider in this paper the sublinear expectation $\ooo,$ which is a functional   defined on a linear space $\mathcal{H}$ and for each $X, Y \in \mathcal{H}$ the following properties are satisfied:\\
\hspace*{0.5cm}(i)  monotonicity: $\ooo [X]\leq \ooo [Y]$ if $X\leq Y;$ \\
\hspace*{0.5cm}(ii) constant preserving: $\ooo [c]=c$ for all   constants $c\in\mathbb{R};$\\
\hspace*{0.5cm}(iii) sub-additivity: $\ooo[X+Y]\le\ooo [X]+\ooo[Y];$\\
\hspace*{0.5cm}(iv) positive homogeneity: $\ooo [\lambda X]=\lambda \ooo[X]$ for all constants $\lambda\ge 0.$\\
The triple $(\Omega, \mathcal{H}, \ooo)$ is called a sublinear expectation space.

As a standard example, consider a sequence $\xi_1,\xi_2,\dots$ of classical  i.i.d. random variables defined on a probability space
$(\Omega, \mathcal{F}, {\bf P}_\theta)$ with $\theta\in\Theta,$ where $\Theta$ contains more than one element.
Denote by $\Psi_n$ the set of  all  bounded and continuous
functions $\psi(\cdot): \mathbb{R}^{n}\rightarrow \mathbb{R},$  and
with $\vec\xi_n=(\xi_1, ,\ldots, \xi_n)$ introduce: %% notation:
\begin{gather}                                                                    \label{a++}
\mathcal{H}:=\{\psi(\vec\xi_n):\psi(\cdot)\in\Psi_n\}\quad\text{and} \quad
\ooo\psi(\vec\xi_n):=\sup_{\theta\in\Theta}{\bf E}_\theta\psi(\vec\xi_n).
\end{gather}
Then  the triple   $(\Omega, \mathcal{H}, \ooo)$  defined   in (\ref{a++}) is a sublinear  expectation  space.

For other examples see \cite{Pe05}  and  \cite{Pe10}.

If we want to obtain CLT or results  about  Prokhorov  distance under sublinear  expectation  then we need a definition of independence.
Similar to Peng  \cite{Pe07b} and \cite{Pe10},
%%%Chen  \cite{C10},
 we  say that a random  vector $\vec{Y}\in \mathbb{R}^l$  is independent to another random vector $\vec{X}\in\mathbb{R}^m$  if
\begin{gather}                                                          \label{a20}
\psi(\vec{X}, \vec{Y})\in\mathcal{H}
%%%, \quad  \psi(\vec x, \vec{Y})\in\mathcal{H}
\quad\text{and} \quad
\ooo\big[\psi(\vec{X}, \vec{Y})\big]=\ooo\big[   \ooo[\psi(\vec x, \vec{Y})]\big|_{\vec x=\vec{X}}\big]
\end{gather}
for all  bounded and
continuous
%%% measurable
function $\psi(\cdot): \mathbb{R}^{l+m}\rightarrow \mathbb{R}.$
This independence means that the distribution  of $\vec{Y}$  does not  change after any realization $\vec x$ of $\vec{X}.$

A sequence of random variables $X_1, X_2, \ldots$  is said to be i.i.d. with $X_0$, if for each $k=2,3,\ldots$ random  variable $X_{k}$ is independent to the vector $\vec{X}_{k-1}=(X_1, X_2,\ldots, X_{k-1})$ and if $X_{k}$
is identically distributed with $X_{0} $  for all $k\ge1$. The latter means that
\begin{gather}                                                          \label{a21}
\forall  k\ge1\ \  \ooo[\psi(X_{k})]=\ooo[\psi(X_0)]\ \  \text{for  any  continuous   and  bounded   } \psi(\cdot).
\end{gather}

Thus,  for a sequence of i.i.d. with $X_0\in\mathcal{H}$ random variables, under a sublinear expectation $\ooo,$ we have two boundaries for values of expectations:
\begin{gather}                                                          \label{a22}
-\infty<\underline{\mu}=\uuu[X_0]=-\ooo[-X_0]\leq\overline{\mu}=\ooo[X_0]<\infty.
\end{gather}
Here and later on we use notation $\uuu[X]:=-\ooo[-X]$ for all $X \in \mathcal{H}.$
And when $X_0^2\in\mathcal{H}$, instead of variance we have two functions:
\begin{gather}\label{a23}
0\le\underline{\sigma}^{2}(\mu)=\uuu[(X_0-\mu)^{2}]
\leq\overline{\sigma}^{2}(\mu)=\ooo[(X_0-\mu)^{2}]<\infty
,\quad \mu\in[\underline{\mu},\overline{\mu}].
\end{gather}

Later on in the paper we everywhere assume that
\begin{gather}                                                          \label{a24}
X_0,X_0^2\in\mathcal{H}
\quad\text{and} \quad
0<\underline{\sigma}:=\inf_{\mu\in[\underline{\mu},\overline{\mu}]}
\underline{\sigma}(\mu)
\le \overline{\sigma}:=\sup_{\mu\in[\underline{\mu},\overline{\mu}]}
\overline{\sigma}(\mu)<\infty.
\end{gather}

{\bf 1.5.}
The situation with CLT under sublinear expectation is essentially  more difficult, even in the case when
\begin{gather}                                                                  \label{a51}
\underline{\mu}=\overline{\mu}=0, \quad
%% \ \text{and}\ \
0<\underline{\sigma}(0)\leq \overline{\sigma}(0)<\infty, \quad
|X_0|^3\in\mathcal{H}.
\end{gather}
Introduce a class of real functions satisfying the H\"older condition:
\begin{gather}                                                            \label{a52}
L_1(\varphi):=\sup\big\{ |\varphi(x+\delta)-\varphi(x)|/\delta:\delta>0,\quad x\in \mathbb{R} \big\}<\infty.
%%\sup_{x_1\neq x_2}\frac{ |\varphi(x_2)-\varphi(x_1)|}{|x_2-x_1|^p}<\infty.
\end{gather}
 In  \cite{Pe08}  and  \cite{Pe10}  Peng  presented the CLT
%%%central limit theorems
for a sequence of  i.i.d. random variables  under sublinear expectations.
Under assumptions (\ref{a51}) and (\ref{a52}) he proved that
 \begin{gather}                                                            \label{a53}
\ooo\varphi(\overline{W}_n)\longrightarrow G(\varphi,c)
\quad\text{for} \quad
\overline{W}_n= \frac{1}{\overline{\sigma}(0)\sqrt{n}}\sum_{i=1}^{n}X_i,\quad c=\frac{\underline{\sigma}(0)}{\overline{\sigma}(0)},
\end{gather}
where the functional $G(\varphi,c)$ may be described in terms of so called
 viscosity solutions of parabolic partial differential  equations. Only when $0<\underline{\sigma}(0)=\overline{\sigma}(0)$ the functional $G(\varphi,1)= {\bf E}\varphi(Z)$
coincides with the right hand side in (\ref{a3}).

The nature of many difficulties  in  understanding of  CLT  under sublinear expectation
 was essentially clarified in the remarkable work \cite{FP17} by Fang, Peng, Shao and Song.
They suggested to consider normalized sums of the following form
\begin{gather}                                                                  \label{a55}
W_{n,n}:= \frac{1}{\sqrt{n}}\sum_{i=1}^{n}\frac{X_i-\mu_i}{\sigma_i},
\end{gather}
where  $\mu_i$ and $\sigma_i$
are  special functions of $\vec{X}_{i-1}$ which
may depend also on numbers $i, n,\varepsilon$  and on the function $\varphi$ with condition (\ref{a52}).
For special class of sublinear expectations with many additional technical assumptions,
including (\ref{a++}) and (\ref{a24}), they proved in  \cite{FP17} that
\begin{gather}                                                                  \label{a56}
 \ooo\varphi(W_{n,n})  -{\bf E}\varphi(Z)  =L_1(\varphi)O\Big(\frac{\log n}{\sqrt{n}}\Big),
\end{gather}
where $Z$ is again the standard Gaussian random variable.
We emphasize that the main advantage of (\ref{a56}) in comparison with (\ref{a53}) consists in the fact that
 the limit ${\bf E}\varphi(Z)$ in (\ref{a56}) is  exactly the same as in the classical CLT.

We are  going to continue  investigations  from \cite{FP17}, but under general definition of sublinear expectations.
And we will consider below expectations of non-continuous indicators instead of smooth functions $\varphi$ with condition (\ref{a52}).

{\bf 1.6.}
In the present paper we consider a sequence of random variables $X_1, X_2, \ldots$  in a sublinear expectation space $(\Omega, \mathcal{H}, \ooo )$, which are supposed to be  i.i.d. with $X_0$, satisfying (\ref{a21})
and (\ref{a24}).
Our main aim  is to generalize  the classical   Prokhorov distances into   subliner expectation space  and obtain  convergence rates for Prokhorov distances in CLT and functional CLT by appropriate generalizations of Lindeberg  method.
For the special broken line $S_n=S_n(\cdot)$ we are going to obtain estimates for the values of the form
\begin{gather}                                                                    \label{a25}
\overline{{\bf P}}(S_n \in D)- {\bf P}(B \in D^{\varepsilon}),   \quad     D \subset \mathbb{C}[0,1],\quad\varepsilon>0,
\end{gather}
where $B=B(\cdot)$ is a Wiener process and
%%where $D^{\varepsilon}$ is defined in (\ref{a16}) and
\begin{gather}                                                                    \label{a26}
\overline{{\bf P}}(S_n \in D):=\inf\{\ooo\psi(\vec{X}_n):
\psi(\vec{X}_n)\ge {\bf I}(S_n \in D), \psi(\vec{X}_n)\in\mathcal{H}\}\le1
\end{gather}
since ${\bf I}(S_n \in D)\le1\in\mathcal{H}$.
Underline that  $\psi(\vec{X}_k)\in\mathcal{H}$  for all   bounded and
continuous functions $\psi(\cdot): \mathbb{R}^{k}\rightarrow \mathbb{R};$ this fact follows from the assumption (\ref{a20}).

The main difference with the classical case is that
 instead of (\ref{a11}) the random broken line $S_n=S_n(\cdot)$ has the following special values at key points:
\begin{gather}                                                                  \label{a27}
S_n(\frac{k}{n})=W_{k,n}:= \frac{1}{\sqrt{n}}\sum_{i=1}^{k}\frac{X_i-\mu_i}{\sigma_i}, \quad i, k=0,1,2,\ldots,n,
\end{gather}
where  $S_n(0)=W_{0,n}=0$ and
%%\begin{gather}                                                                    \label{a28}
$\mu_i=\mu_{i,n,D,\varepsilon}(\vec{X}_{i-1})$ and
%%\quad\text{and} \quad
$\sigma_i=\sigma_{i,n,D,\varepsilon}(\vec{X}_{i-1})$ %%\  \  i=1,2,\ldots,n,
%%\end{gather}
are  special continuous functions of $\vec{X}_{i-1}$ which
may depend also on numbers $i, n,\varepsilon$  and on the set $D$ from (\ref{a25}). In addition
\begin{gather}                                                                    \label{a29}
\underline{\mu}\le\mu_i\le\overline{\mu}
\quad\text{and} \quad
0<\underline{\sigma}\le\underline{\sigma}(\mu_i)\le\sigma_i\le\overline{\sigma}(\mu_i)\le\overline{\sigma}
, \quad i=1,2,\ldots,n.
\end{gather}

In Section 2 we introduce Prokhorov distances in sublinear expectation spaces and present estimates for them as  for one-dimensional CLT, so  for functional CLT. Underline, that our estimates from Theorem \ref{thm2} in functional CLT is, in some sense, better than unimprovable Borovkov's estimate (\ref{a41}). Indeed, for the classical Prokhorov distance under standard assumption (\ref{a1}) our Theorem \ref{thm2} yields, for $p\in[2,3]$, the following inequality
\begin{gather}                                                                                                              \label{a42}
\Pi(S_n, B)\le 4.7\frac{\big({\bf E}\min\big\{|\xi_1|^p,\sqrt n\xi_1^2\big\}\big)^{1/(p+1)}
}{n^{(p-2)/(2p+2)}}
\le 4.7\frac{\big({\bf E}|\xi_1|^p\big)^{1/(p+1)}
}{n^{(p-2)/(2p+2)}}.
%%,\quad p\in[2,3].
\end{gather}
Thus, we obtained Borovkov's estimate (\ref{a41}) with constant $C_1(p)=4.7$ for $p\in[2,3]$.
As far as we know it is the first case when an estimate in the classical functional CLT is obtained with an explicit numerical constant (about previous attempt see Remark \ref{R41} below). Moreover, the first inequality in (\ref{a42}) is even better, because it uses truncated moment and  does not require the existence of moments, but variance.

{\bf 1.7.}
The rest of the paper is organized as follows.  The main results  about Prokhorov  distances defined on   sublinear  expectation spaces  are presented  in Section~2.  In  Sections~3 and~4 we  prove our main Theorem~\ref{thm2} and its corollaries. In section 5  we  prove simpler Theorem~\ref{thm1CLT}
with estimates in one-dimensional~CLT.

\section{ Main results}
{\bf 2.1. Preliminary agreements and notations.}
Later on for any function  $\phi(\cdot): \mathbb{R}\rightarrow \mathbb{R}$ we  use notation
\begin{gather*}      %                                                              \label{b1}
\ooo[\phi(X_0)]:=\inf\big\{\ooo[X]:
X\ge \phi(X_0) , X\in\mathcal{H}\big\}\le\infty,
\end{gather*}
where we assume that $\inf\emptyset=\infty$. Now for
$\mu\in[\underline{\mu},\overline{\mu}]$, $p\ge2$ and $C>0$ define
\begin{gather}                                                          \label{b2}
\xi_0(\mu):=(X_0-\mu)/\underline{\sigma}(\mu)\quad\text{and} \quad
\gamma_p(\mu,C):=
\ooo\big[|\xi_0(\mu)|^p\wedge C^{p-2}\xi_0^2(\mu)\big],
%\quad p\ge2,
\end{gather}
where $a\wedge b$ means $\min\{a,b\}$. Introduce
\begin{gather}                                                         \label{b3}
\gamma_p(C):=\sup_{\mu\in[\underline{\mu},\overline{\mu}]}
\gamma_p(\mu,C)\le C^{p-2} {\overline{\sigma}^2}/{\underline{\sigma}^p}<\infty ,\quad C>0,\quad p\ge2.
\end{gather}
Define also
\begin{gather*}           %                                               \label{b4}
\sup_{C>0}
\gamma_p(C)\le\gamma_p:= \sup_{\mu\in[\underline{\mu},\overline{\mu}]}
\ooo\big[|\xi_0(\mu)|^p\big]\le\infty ,\quad p\ge2.
\end{gather*}

\begin{rem}                                                        \label{R21}
We treat the distribution of $X_0$ as a specially chosen
%approximately the best
one in the class of distributions, which satisfy (\ref{a21}) and (\ref{a24}),
and which approximately minimize  in this class the value $\gamma_p$ or $\gamma_p(C)$ for some $p$ and $C$.
In particular, we do not assume, for example,  that $\ooo X_i^2<\infty$ for some $i\neq0$.
\end{rem}

\noindent{\bf 2.2. Estimates in functional CLT.}
For each $i\ge1$ introduce into consideration the set $\mathcal{M}(i)$ of all pares
 $(\mu_i, \sigma_i)$ which are
continuous functions of $\vec{X}_{i-1}$ and satisfy conditions (\ref{a29}).
And denote by $\mathcal{M}_n$ the set $\big\{(\mu_i, \sigma_i)\in\mathcal{M}(i),\ i=1,\dots, n \big\}$
of all sequences of length $n$ of such pares, which are essentially used in construction (\ref{a27}) of random broken line $S_n$.
\begin{thm}\label{thm2}
  Under assumptions (\ref{a24})  for  each set $D \subset \mathbb{C}[0,1]$ and any  $\varepsilon>0$  there exist   values
$\big\{(\mu_i, \sigma_i)\big\}\in\mathcal{M}_n$ which may depend on $D, \varepsilon, n$  and such that
\begin{gather}                                                                                                           \label{b11}
\overline{{\bf P}}(S_n \in D)- {\bf P}(B \in D^{\varepsilon})
\leq\tilde \Pi_n(\varepsilon,p):=
\frac{C_2(p)\gamma_p(\varepsilon\sqrt{n})}{n^{(p-2)/2}\varepsilon^{p}}
\leq\frac{C_2(p)\gamma_p}{n^{(p-2)/2}\varepsilon^{p}},
%%\quad p\in[2,3],
\end{gather}
 where $p\in[2,3]$ and $C_2(p):=\min\big\{184,(4.7)^{p+1}\big\}$.
\end{thm}

In other words, we proved in (\ref{b11}) that for $p\in[2,3]$
\begin{gather*}       %                                                                                                    \label{b13}
\overline\Delta_n(D,\varepsilon):=\inf_{\{(\mu_i, \sigma_i)\}\in\mathcal{M}_n}
\overline{{\bf P}}(S_n \in D)- {\bf P}(B \in D^{\varepsilon})
\leq\tilde \Pi_n(\varepsilon,p).
\end{gather*}
Introduce notation:
\begin{gather}                                                                                                       \label{b15}
%%\overline{\Pi}(S_n, B, \varepsilon)=\sup_{D\subset \mathbb{C}[0,1]}\overline\Delta_n(D,\varepsilon),\  \
\overline{\Pi}(S_n, B)=\inf\big\{\varepsilon>0:\sup_{D\subset \mathbb{C}[0,1]}\overline\Delta_n(D,\varepsilon)
% \overline{\Pi}(S_n, B, \varepsilon)
\leq \varepsilon\big\}.
\end{gather}
 Later on we consider the value $\overline{\Pi}(S_n, B)$ as Prokhorov distance between distributions
in $\mathbb{C}[0,1]$ of $S_n$ and  $B$ defined on a sublinear expectation space.

\begin{cor}\label{cor1+}
Under conditions of Theorem \ref{thm2} %, we further have,
\begin{gather}                                                                                              \label{b17}
\overline{\Pi}(S_n, B)
%=\inf\big\{\varepsilon>0, \overline{\Pi}(S_n, B, \varepsilon)\leq \varepsilon\big\}
\leq\tilde{\Pi}_{n,p}:= \frac{C_3\big(\gamma_p(\sqrt n)\big)^{1/(p+1)}}{n^{(p-2)/(2p+2)}}
\le\frac{C_3\gamma^{1/(p+1)}_p}{n^{(p-2)/(2p+2)}},
\end{gather}
where $p\in[2,3]$ and $C_3=4.7$.
\end{cor}

As it was mentioned in introduction, in classical case our estimate (\ref{b17}) for our generalization of
Prokhorov distance implies known unimprovable estimate~(\ref{a41}) of Borovkov.
Note also that our Prokhorov distance satisfies the following analog of an classical property (\ref{a19}).

\begin{cor}\label{cor2+}
Under assumptions of Theorem \ref{thm2}, for each functional $f(\cdot)$ which
% is  bounded  and  Lipschitz  continuous, i.e., $f(\cdot)$
 satisfies conditions (\ref{a17}) and (\ref{a18}) we have
\begin{gather}                                                                                              \label{b19}
\inf_{\{(\mu_i, \sigma_i)\}\in\mathcal{M}_n}\overline{{\bf P}}(f(S_n)\leq a)- {\bf P}(f(B)\leq a)
%%\leq\frac{C_4\gamma_3^{1/4}(KL)^{3/4}}{n^{1/8}}
\leq (KL+1)
\tilde{\Pi}_{n,p}.
%%\overline{\Pi}(S_n, B).
\end{gather}
\end{cor}

\noindent{\bf 2.3. Estimates in one-dimensional  CLT.}
For every $i\ge1$ introduce  the set $\mathcal{A}(i)$ of all pares
 $(\mu_i, \sigma_i)$ which are (one-dimensional)
continuous functions of $W_{i-1,n}$ and satisfy conditions (\ref{a29}).
Denote by $\mathcal{A}_n$ the set $\big\{(\mu_i, \sigma_i)\in\mathcal{A}(i),\ i=1,\dots, n \big\}$
of all sequences of length $n$ of such pares.
\begin{thm}                                                                             \label{thm1CLT}
Under assumptions (\ref{a24})  for  each set $A\subset \mathbb{R}$ and any  $\varepsilon>0$  there exist   values
$\{(\mu_i, \sigma_i)\}\in\mathcal{A}_n$ which may depend on $A, \varepsilon, n$  and such that
\begin{gather}                                                                                                           \label{b21}
\overline{{\bf P}}(W_{n,n} \in A)- {\bf P}(Z \in A^{\varepsilon})
\leq
%%\tilde \pi_n(\varepsilon,p):=
\frac{C_4\gamma_p(\varepsilon\sqrt{n})}{n^{(p-2)/2}\varepsilon^{p}}
\leq\frac{C_4\gamma_p}{n^{(p-2)/2}\varepsilon^{p}},
\quad p\in[2,3],
\end{gather}
 where $Z\sim N(0,1)$
%%,$ $p\in[2,3]$
and $C_4=42<(3.5)^{p+1}$.
In addition, when
%% $0<\varepsilon\le1$ and
$\gamma_3<\infty$
we also have
\begin{gather}                                                           \label{b22}
\overline{{\bf P}}(W_{n,n}\in A)-{\bf P}(Z\in A^{\varepsilon})\leq \frac{C_5\gamma_3}{\varepsilon\sqrt{n}}
\quad\text{with} \quad C_5=12.
\end{gather}
%%where $C_5=29.$
\end{thm}

It is equivalent to say that we obtained in (\ref{b21}) and (\ref{b22}) estimates for the value
\begin{gather*}       %                                                                                                    \label{b23}
\overline\delta_n(A,\varepsilon):=\inf_{\{(\mu_i, \sigma_i)\}\in\mathcal{A}_n}
\overline{{\bf P}}(W_{n,n} \in A)- {\bf P}(Z \in A^{\varepsilon}).
%%\leq\tilde \Pi_n(\varepsilon,p).
\end{gather*}
Introduce notation:
\begin{gather}                                                                                                       \label{b25}
%%\overline{\pi}(W_n, Z, \varepsilon)=\sup_{A\subset \mathbb{R}}\overline\delta_n(A,\varepsilon),\quad
\overline{\pi}(W_{n,n}, Z)=\inf\big\{\varepsilon>0: % \overline{\pi}(W_n, Z, \varepsilon)
\sup_{A\subset \mathbb{R}}\overline\delta_n(A,\varepsilon)\leq \varepsilon\big\}.
\end{gather}
 Later on we treat
%%consider the value
$\overline{\pi}(W_{n,n}, Z)$ as Prokhorov distance between distributions
in $\mathbb{R}$ of $W_{n,n}$ and  $Z$ defined on a sublinear expectation space.

We have immediately from (\ref{b22}) and (\ref{b25}) that
\begin{gather}                                                                                                        \label{b28}
\overline{\pi}(W_{n,n}, Z)
%%=\inf\big\{\varepsilon>0, \overline{\pi}(W_{n,n}, Z, \varepsilon)\leq \varepsilon\big\}
\leq
%\tilde\pi_n:=
{\sqrt{C_5\gamma_3}}\big/{n^{1/4}}.
\end{gather}

\begin{rem}
For classical Prokhorov distance ${\pi}(W_{n}, Z)$ in classical CLT for i.i.d.  random  variables under assumptions (\ref{a1})
Yurinsky \cite{V75}
%(1975)
obtained the following estimate
\begin{gather*}                     %                                                                                   \label{b30}
{\pi}(W_{n}, Z)
\leq{C_6{\bf E}|\xi_1|^3}\big/{\sqrt n},\qquad\xi_1:=(X_{1}-\mu)/{\sigma},
\end{gather*}
where $C_6$ is an absolute   constant. Unfortunately, we do not see a method to obtain such sharp estimate for
$\overline{\pi}(W_{n,n}, Z)$, instead of (\ref{b28}).
\end{rem}

Underline also that in Theorems \ref{thm2} and \ref{thm1CLT} we never supposed that sets~$D$ or~$A$ be measurable.

The rest of the paper is devoted to proofs of results presented above.

\section{Main lemmas }

{\bf 3.1. Special functions.}
Let  $X_0$  be  a random  variable  under  sublinear  expectation  which  satisfies  condition  (\ref{a24}).  Then  for  all twice  continuously  differentiable function  $h(\cdot)$ and for  all  non-random  $\varkappa >0$ we  may  introduce  the  following  functions:
\begin{gather}  \label{c1}
\alpha_{\varkappa }(x|h):=\frac{\overline{u}+\underline{\mu}}{2}+\frac{\overline{u}-\underline{\mu}}{2}\cdot \frac{h^{'}(x)}{\varkappa },\  \   \text{ if }\ \ |h^{'}(x)|\le  \varkappa ;
\end{gather}
and when  $|h^{'}(x)|> \varkappa $  we  put
\begin{gather}                                                                                               \label{c2}
\alpha_{\varkappa }(x|h):=
\begin{cases}
\overline{\mu},    \ \text{if}\ \ h^{'}(x)> \varkappa ,\\
\underline{\mu}, \ \text{if}\ \  h^{'}(x)<-\varkappa .
\end{cases}
\end{gather}
It  is  easy  to  see  that
$\alpha_{\varkappa }(x|h)$  is  continuous  for  all  $x\in \mathbb{R}$, since it is a continuous function of $h^{'}(x)$.

For  all  non-random  $\kappa >0$  introduce
\begin{gather}\label{c3}
\beta^2_{\varkappa ,\kappa }(x|h):=\frac{\overline{\sigma}^2(\alpha_{\varkappa }(x|h))+\underline{\sigma}^2(\alpha_{\varkappa }(x|h))}{2}+
\frac{\overline{\sigma}^2(\alpha_{\varkappa }(x|h))-\underline{\sigma}^2(\alpha_{\varkappa }(x|h))}{2}\cdot \frac{h^{''}(x)}{\kappa },
%%\quad   \text{  if }  |h^{''}(x)|\le \kappa ;
\end{gather}
when $|h^{''}(x)|\le \kappa $;
and if $|h^{''}(x)|>\kappa $ we put
\begin{gather}\label{c4}
\beta^2_{\varkappa ,\kappa }(x|h):=
\begin{cases}
\overline{\sigma}^2(\alpha_{\varkappa }(x|h)),    \ \text{if}\ \ h^{''}(x)>\kappa ,\\
\underline{\sigma}^2(\alpha_{\varkappa }(x|h)), \ \text{if}\ \  h^{''}(x)<-\kappa .
\end{cases}
\end{gather}
It is  easy to  see  that $\beta_{\varkappa ,\kappa }(x|h)$  is  also  continuous  for  all   $x\in \mathbb{R}$,
because it is a continuous function of $h^{''}(x)$ and $h^{'}(x).$
In addition, for  all   $x\in \mathbb{R}$
\begin{gather}                                                                    \label{c5+}
\underline{\mu}\le\alpha_{\varkappa }(x|h)\le\overline{\mu}
\quad\text{and} \quad
0<\underline{\sigma}\le\underline{\sigma}(\alpha_{\varkappa }(x|h))\le\beta_{\varkappa ,\kappa }(x|h)\le\overline{\sigma}(\alpha_{\varkappa }(x|h))\le\overline{\sigma}.
\end{gather}

For  any   $\varkappa , \kappa >0$  and  $x\in  \mathbb{R}$  define
\begin{gather}\label{c5}
\xi(x)=(X_0-\alpha_{\varkappa }(x|h))/\beta_{\varkappa ,\kappa }(x|h);
%\frac{X_0-\alpha_{\varkappa }(x|h)}{\beta_{\varkappa ,\kappa }(x|h)},
\quad\text{so that} \quad
|\xi(x)|\le|\xi_0(\alpha_{\varkappa }(x|h))|.
\end{gather}
\begin{lem}\label{LL-1}
If  $h(\cdot)$  is  a twice  continuously  differentiable function, then  for  all $\varkappa , \kappa >0$  and  $x\in \mathbb{R}$
\begin{eqnarray}\label{c6}
\ooo h(x+\frac{\xi(x)}{\sqrt{n}})&\le& {\bf E}h(x+\frac{Z}{\sqrt{n}})
+\frac{\varkappa (\overline{\mu}-\underline{\mu})}{2\sqrt{n}\underline{\sigma}}
+\frac{\kappa }{2n}(\frac{\overline{\sigma}^2}{\underline{\sigma}^2}-1)\nonumber\\
&+&\ooo\delta(x, \frac{\xi(x)}{\sqrt{n}}| h)-{\bf E}\delta(x, \frac{Z}{\sqrt{n}}|h ),
\end{eqnarray}
where classical  random variable  $Z$  has  standard  normal  distribution  and
\begin{gather} \label{c7}
\delta(x,y|h)=h(x+y)-h(x)-yh^{'}(x)-{y^2}h^{''}(x)/2.
\end{gather}
\end{lem}

\begin{rem}\label{R31}
If  $X_0$  is a  classical random variable  we have
$$\xi(x)=\frac{X_0-{\bf E}X_0}{\sqrt{{\bf Var} X_0}},\  \  {\bf E} \xi(x)={\bf E} Z=0, \ \ {\bf E} \xi^2(x)={\bf E} Z^2=1.$$
Hence,  instead  of  (\ref{c6}) we  have
\begin{gather}  \label{c11}
{\bf E}h(x+\frac{\xi(x)}{\sqrt{n}})={\bf E}h(x+\frac{Z}{\sqrt{n}})+{\bf E} \delta(x, \frac{\xi(x)}{\sqrt{n}}|h)-{\bf E}\delta(x, \frac{Z}{\sqrt{n}}|h).
\end{gather}
Equality (\ref{c11})  plays  an important  role  in the  Lindeberg method  to  prove CLT.  By  this  reason  we can consider  our  proof,  based on  Lemma  \ref{LL-1}, as a natural generalization  of  the Lindeberg  method  for the  case of  sublinear  expectations.
\end{rem}
\begin{rem}\label{R32}
We  may define  functions  $\alpha_0(x|h)$  and  $\beta_{0,0}(x|h)$  by  formulas (\ref{c2})  and  (\ref{c4})  with  $\varkappa =\kappa =0.$  In this  case  we  may  expect  that for  all  $x\in \mathbb{R},$
\begin{gather}                                                                                                            \label{c12}
\ooo\big[  h^{'}(x)\xi(x)\big]=0\ \  \text{and }\  \  \ooo\big[h^{''}(x)\xi^2(x)   \big]=h^{''}(x).
\end{gather}
But
%the  main  difficulty  in  proving  of  (\ref{c12})  is  that
 functions  $\alpha_0(x|h)$  and $\beta_{0,0}(x|h)$  are  not  continuous.  By  this  reason, in this  paper  we use  more complicated  continuous  functions $\alpha_{\varkappa }(x|h)$  and $\beta_{\varkappa ,\kappa }(x|h)$ instead  of the simpler  functions   $\alpha_0(x|h)$  and $\beta_{0,0}(x|h).$
\end{rem}

\noindent{\bf 3.2. Proof of Lemma \ref{LL-1}.}
 For  a  fixed $x$ take  Taylor  expansions  of
%%$h(x+\frac{\xi(x)}{\sqrt{n}})$  and
$h(x+\frac{Z}{\sqrt{n}})$  at $x$.
%%respectively,  t
Then
\begin{gather}\nonumber
h(x+\frac{Z}{\sqrt{n}})=h(x)+\frac{Z}{\sqrt{n}}h^{'}(x)+\frac{Z^2}{2n}h^{''}(x)+\delta(x, \frac{Z}{\sqrt{n}}|h),
\\                                                                                                           \label{c41}
{\bf E} h(x+\frac{Z}{\sqrt{n}})=h(x)+\frac{1}{2n}h^{''}(x)+{\bf E}\delta(x, \frac{Z}{\sqrt{n}}|h).
\end{gather}
Similarly,
$$h(x+\frac{\xi(x)}{\sqrt{n}})=h(x)+\frac{\xi(x)}{\sqrt{n}}h^{'}(x)+\frac{\xi^2(x)}{2n}h^{''}(x)+\delta(x, \frac{\xi(x)}{\sqrt{n}}|h),$$
and by the properties (ii) -- (iv)  of  sublinear  expectation we  have
\begin{gather}                                                                                    \label{c42}
\ooo\big[h(x+\frac{\xi(x)}{\sqrt{n}})\big]
\le h(x)+\frac{\ooo\big[ \xi(x)h^{'}(x) \big]}{\sqrt{n}}+\frac{\ooo\big[ \xi^2(x)h^{''}(x)  \big]}{2n}
+\ooo \delta(x, \frac{\xi(x)}{\sqrt{n}}|h).
\end{gather}

We are going to estimate values $\ooo\big[ \xi(x)h^{'}(x)  \big]$ and $\ooo\big[ \xi^2(x)h^{''}(x)  \big]$.
To simplify notations, we will write in the proof $\alpha(x)$ and $\beta(x)$ instead of
$\alpha_{\varkappa }(x|h)$ and $\beta_{\varkappa ,\kappa }(x|h)$.
We will use the following property
\begin{gather}                                                                                               \label{c43}
\ooo\big[c(X-a)\big]=
\begin{cases}
c(\ooo[X]-a)=|c|(\ooo[X]-a),    \ \text{if}\ \ c\ge 0,\\
|c|(a+\ooo[-X])=|c|(a-\uuu[X]), \ \text{if}\ \  c\le0,
\end{cases}
\end{gather}
which takes place for all constants $a,c$   and any variable $X\in\mathcal{H}$ in  a sublinear expectation space $(\Omega, \mathcal{H}, \ooo )$.

First, if  $h^{'}(x)\ge 0$ then, by (\ref{c43})
$$\ooo\big[ (X_0-\alpha(x))h^{'}(x) \big]=|h^{'}(x)|(\overline{\mu}-\alpha(x)).$$
But  condition  $\alpha(x)\neq\overline{\mu}$  when  $h^{'}(x)\ge 0$ implies that $|h^{'}(x)|\le \varkappa $
  and  $\alpha(x)\ge ({\overline{\mu}+\underline{\mu}})/{2}.$
Hence
$$\ooo\big[ (X_0-\alpha(x))h^{'}(x) \big]\le\varkappa  (\overline{\mu}-\alpha(x))
\le{\varkappa (\overline{\mu}-\underline{\mu})}/2,\ \ \text{ when } h^{'}(x)\ge 0.$$
Similarly,  when  $h^{'}(x)< 0,$ using again (\ref{c43}), we obtain that
$$\ooo\big[ (X_0-\alpha(x))h^{'}(x) \big]=|h^{'}(x)|(\alpha(x)-\uuu[X_0])
=|h^{'}(x)|(\alpha(x)-\underline{\mu}).$$
But  by  (\ref{c1}) and   (\ref{c2}), condition  $\alpha(x)\neq\underline{\mu}$  when  $h^{'}(x)<0$  implies  that
$|h^{'}(x)|\le \varkappa $
  and  $\alpha(x)\le ({\overline{\mu}+\underline{\mu}})/{2}.$
Hence
$$\ooo\big[ (X_0-\alpha(x))h^{'}(x) \big]\le\varkappa (\alpha(x)-\underline{\mu})
\le {\varkappa (\overline{\mu}-\underline{\mu})}/{2},\ \ \text{ when } h^{'}(x)< 0.$$
Thus,  for  all  values of $h^{'}(x)$ we have
\begin{gather}                                                                                          \label{c44}
\ooo\big[ \xi(x)h^{'}(x) \big]=\frac{1}{\beta(x)}\ooo\big[ (X_0-\alpha(x))h^{'}(x)  \big] \le \frac{\varkappa (\overline{\mu}-\underline{\mu})}{2\underline{\sigma}}.
\end{gather}

Next,
%%we are  going to  estimate $\overline{E}\big[ h^{''}(x)\xi^2(x) \big].$
if  $h^{''}(x)\ge 0$
then, by (\ref{c5}),
$$\ooo\big[ h^{''}(x)\xi^2(x) \big]
=h^{''}(x)\ooo\big[ \frac{(X_0-\alpha(x))^2}{\beta^2(x)}   \big]
=h^{''}(x)\frac{\overline{\sigma}^2(\alpha(x))}{\beta^2(x)}
%%%\le h^{''}(x)\frac{\overline{\sigma}^2}{\beta^2(x)}
$$
%%since $\beta^2(x)\le \overline{\sigma}^2$ for all $x$ by (c5+) and (a24).
But   by (\ref{c4}) and (\ref{c3}), condition  $\beta^2(x)\neq \overline{\sigma}^2(\alpha(x))$  when $h^{''}(x)\ge 0$  implies that
$h^{''}(x)\le \kappa $.
Hence
\begin{gather}                                                                     \label{c45+}
\ooo\big[h^{''}(x)\xi^2(x)\big]-h^{''}(x)= h^{''}(x)(\frac{\overline{\sigma}^2(\alpha(x))}{\beta^2(x)}-1)
%%\frac{\kappa }{2}(\frac{\overline{\sigma}^2(\alpha(x))}{\underline{\sigma}^2(\alpha(x))}-1)
\le
\kappa (\frac{\overline{\sigma}^2}{\underline{\sigma}^2}-1),
\  \   \text{if }\  \  h^{''}(x)\ge 0.
\end{gather}

At last, when  $h^{''}(x)<0,$  we     have from    (\ref{c5})  and  (\ref{c43})   that
\begin{eqnarray*}
\ooo\big[ h^{''}(x)\xi^2(x) \big]-h^{''}(x)
=|h^{''}(x)|(- \uuu\big[ \xi^2(x) \big]+1)
\\ \nonumber
=|h^{''}(x)|(1- \uuu\big[
 \frac{(X_0-\alpha(x))^2}{\beta^2(x)}   \big])
=|h^{''}(x)|(1-\frac{\underline{\sigma}^2(\alpha(x))}{\beta^2(x)}).
%%%\le |h^{''}(x)|(1-\frac{\underline{\sigma}^2}{\beta^2(x)}).
\end{eqnarray*}
%%%since $\beta^2(x)\ge \underline{\sigma}^2$ for all $x$ by (c5+) and (a24).
By (\ref{c4}) and (\ref{c3}), condition  $\beta^2(x)\neq\underline{\sigma}^2(\alpha(x))$  when $h^{''}(x)< 0$  implies that
$h^{''}(x)\le \kappa $.
%%%  and  $\beta^2(x)\le  \overline{\sigma}^2(\alpha(x))$
 Hence
\begin{gather}                                                                                                    \label{c45-}
\ooo\big[h^{''}(x)\xi^2(x)\big]-h^{''}(x)\le \kappa (1-\frac{\underline{\sigma}^2(\alpha(x))}{\beta^2(x)})
\le
\kappa (1-\frac{\underline{\sigma}^2}{\overline{\sigma}^2}),
  \text{  when  }  h^{''}(x)< 0.
\end{gather}
%%since $\beta^2(x)\ge \overline{\sigma}^2(\alpha(x))$  for  all   $x.$

Note that in (\ref{c45+}) and (\ref{c45-}) we also used (\ref{c5+}).
But $1-{\underline{\sigma}^2}/{\overline{\sigma}^2}\le{\overline{\sigma}^2}/{\underline{\sigma}^2}-1$.
Thus   for  all  values of $h^{''}(x)$
\begin{gather}\label{c45}
\ooo\big[h^{''}(x)\xi^2(x)\big]-h^{''}(x)
\le \kappa ({\overline{\sigma}^2}/{\underline{\sigma}^2}-1).
\end{gather}
Substituting now (\ref{c44}) and (\ref{c45}) into (\ref{c42}), we obtain that
\begin{eqnarray}\label{c46}
\ooo h(x+\frac{\xi(x)}{\sqrt{n}})\le h(x)
+\frac{\varkappa (\overline{\mu}-\underline{\mu})}{2\sqrt{n}\underline{\sigma}}
+\frac{\kappa }{2n}(\frac{\overline{\sigma}^2}{\underline{\sigma}^2}-1)
+\ooo\delta(x, \frac{\xi(x)}{\sqrt{n}}|h).
\end{eqnarray}
So, the desired inequality
  (\ref{c6}) follows from     (\ref{c41})    and     (\ref{c46}).

\medskip\noindent{\bf 3.3. Smoothing.}
 \ Introduce  functions
\begin{gather}                                                                                                  \label{g1}
g_r(x)=\frac{2((r-|x|)^{+})^2-((r-2|x|)^{+})^2}{r^3}, \   \  x\in \mathbb{R}, \  \  r>0.
\end{gather}
\begin{lem}                                                                                                                \label{LL-3}
For   each  $r>0$  function  $g_r(x)$  is  a continuous density  on  $\mathbb{R}$  with  a derivative  $g^{'}_r(x)$  which  is absolutely continuous on  $\mathbb{R}.$  It's  second derivative $g^{''}_r(x)$  is  also  continuous  for  all   $x\in \mathbb{R},$   but  $x=\pm r$  and  $x=\pm r/2.$  Moreover,
\begin{gather}                                                                                                \label{g2}
%%\overline{\varepsilon}(y|g_r)=
\overline{\varepsilon}_{r}(y):=\frac12\int |\varepsilon_{r}(x,y)|dx
\le 16 \min\big\{\frac{ |y|^3}{(2r)^3}, \frac{ y^2}{(2r)^2}\big\}
\le 16 \min\big\{\frac{ |y|^p}{(2r)^p}, \frac{ y^2}{(2r)^2}\big\}
,
\end{gather}
for  all  $r>0$, $p\in[2,3]$  and  $y\in \mathbb{R}$, where
\begin{gather}                                                                                          \label{g3}
\varepsilon_{r}(x,y):=g_r(x-y)-g_r(x)+yg^{'}_r(x)-{y^2g^{''}_r(x)}/{2}.
\end{gather}

\end{lem}
{\bf  Proof. }  If $r=1/2$  then  the assertion follows  from  the  proof   of  Lemma 6.2  in  Sakhanenko (2000).
The  case  of  arbitrary  $r>0$   follows  from  the evident representations:
\begin{gather*}%%\label{d6}
g_r(x)=\frac{1}{2r}g_{1/2}(\frac{x}{2r}), \quad
% g^{'}_r(x)=\frac{1}{(2r)^2}g^{'}_{1/2}(\frac{x}{2r}),
%\quad   g^{''}_r(x)=\frac{1}{(2r)^3}g^{''}_{1/2}(\frac{x}{2r}),
%\\ %\end{gather}\begin{gather}\label{d7}
\varepsilon_{r}(x,y)=\frac{1}{2r}\varepsilon_{1/2}(\frac{x}{2r}, \frac{y}{2r}),\  \ \quad
\overline{\varepsilon}_{r}(y)=\overline{\varepsilon}_{1/2}( \frac{y}{2r}).
\end{gather*}
\hfill\fbox

Suppose now that for each $\vec x\in R^m$ we are given a measurable set $A(\vec x)\subset \mathbb{R}$.
\begin{lem}                                                                                                                \label{LL-2}
Let random vector $\vec\eta$ be independent from random variable $\nu$ with density $g_1(\cdot)$.
Then the following function
\begin{gather*}      %                                                                                    \label{g5}
h(x):={\bf P}\big(x+r\nu\in A(\vec\eta)\big),\quad r>0,
\end{gather*}
 is twice continuously differentiable for all $x\in \mathbb{R}$ and
%%, moreover, for each $b\neq 0$ and $y\in \mathbb{R}$
\begin{gather}                                                                                                \label{g6}
\overline{\delta}(y|h):=\sup_{x\in \mathbb{R}}\big|\delta(x,y|h)\big|\le \overline{\varepsilon}_r(y),
\end{gather}
for each $r> 0$ and $y\in \mathbb{R}$, where the function $\delta(x,y|h)$ was defined in (\ref{c7}).
\end{lem}

{\bf  Proof. }
Remind that for each non-random $\vec{s}\in \mathbb{R}^m$
\begin{gather*} %%                                                                                \label{c25}
 {\bf P}\big(x+r\nu\in A(\vec{s})\big)=\int_{A(\vec{s})}g_r(z-x)dz,
\end{gather*}
where  $g_r(\cdot)$  is the  density  of  $r\nu.$
Now for any $\vec\eta$ independent from  %%random variable
$\nu$ we have
\begin{gather}                                                                                \label{g7}
h(x)= {\bf P}\big(x+r\nu\in A(\vec\eta)\big)={\bf E}\int_{A(\vec\eta)}g_r(z-x)dz
\end{gather}
From definitions (\ref{c7}), (\ref{g3}) and (\ref{g7}) it is not difficult to obtain that
%%of  $\delta(x,y|g)$  in (\ref{c7})  and  (\ref{c22}) we have
\begin{gather}                                                                                \label{g11}
\delta(x, y|h)={\bf E}\int_{A(\vec\eta)}
\varepsilon_r(z-x, y)dz.
\end{gather}

Note
%It follows from (c20)
that $$\int\varepsilon_r(x, y)dx=0=\int\varepsilon_r^+(x, y)dx-\int\varepsilon_r^-(x,y)dx.$$
Hence, we have from (\ref{g11}) that
\begin{gather*}                                                                                %%\label{c28}
\delta(x,y|h)\le \int\varepsilon_r^+(z-x, y)dz
=\int\varepsilon_r^+(x,y)dx=\frac12\int|\varepsilon_r(x, y)dx|dx,
\\     %%                                                                           \label{c29}
-\delta(x,y|h)\le \int\varepsilon_r^-(z-x, y)dz
=\int\varepsilon_r^-(x,y)dx=\frac12\int|\varepsilon_r(x, y)dx|dx.
\end{gather*}
These two inequalities imply (\ref{g6}).
\hfill\fbox

\begin{lem}                                                                                                                \label{LL-4}
Let random vector $\vec\eta$ be independent from random variable $Z\sim N(0,1).$
Then the following function
\begin{gather*}         %                                                                                 \label{g15}
h(x):={\bf P}\big(x+bZ\in A(\vec\eta)\big),\qquad b\neq0,
\end{gather*}
 is twice continuously differentiable for all $x\in \mathbb{R}$ and
%%, moreover, for each $b\neq 0$ and $y\in \mathbb{R}$
\begin{gather}                                                                                                \label{g16}
\overline{\delta}(y|h):=\sup_{x\in \mathbb{R}}\big|\delta(x,y|h)\big|\le 0.4|y/b|^3
%%\frac{2|y|^3}{5|b|^3},
\end{gather}
for each $y\in \mathbb{R}$, where the function $\delta(x,y|h)$ was defined in (\ref{c7}).
\end{lem}

{\bf  Proof. }
Similarly to the proof of Lemma \ref{LL-2} we have
\begin{gather}                                                                                \label{g17}
h(x)= {\bf P}\big(x+bZ\in A(\vec\eta)\big)={\bf E}\int_{A(\vec\eta)}p_b(z-x)dz
\end{gather}
where  $p_b(\cdot)$  is the  density  of  $bZ.$
Now introduce
\begin{gather}                                                                                         \label{g18}
\varepsilon(x,y):=p_b(x-y)-p_b(x)+yp^{'}_b(x)-\frac{y^2p^{''}_b(x)}{2}=\frac{y^3}{2}\int^{1}_0 (1-t)^2p^{'''}_b(x-ty)dt.
\end{gather}
From definitions (\ref{c7}), (\ref{g17}) and (\ref{g18}) it is not difficult to obtain that
\begin{gather*}       %                                                                         \label{g19}
\delta(x, y|h)={\bf E}\int\limits_{A(\vec\eta)}\varepsilon(z-x, y)dz
=\frac{y^3}{2}\int^1_0 (1-t)^2{\bf E}\int\limits_{A(\vec\eta)}p^{'''}_b(z-x-ty)dzdt.
\end{gather*}
So,  we have that
\begin{gather}
\big|\delta(x, y|h)\big|\le \frac{|y|^3}{2}\int^1_0(1-t)^2\int |p^{'''}_b(z-x-ty)|dzdt\nonumber
\\                                                                                       \label{g20}
=\frac{|y|^3}{2}\int^1_0(1-t)^2dt\cdot \int |p^{'''}_b(y)|dy=\frac{|y|^3}{6|b|^3}\int |p^{'''}_1(x)|dx.
\end{gather}

It  is  easy  to  see that
\begin{gather*}   %                                                                                      \label{g21}
p^{'}_1(x)=-xp_1(x), \ \ p^{''}_1(x)=(x^2-1)p_1(x)\   \   \text{and }\  \  p^{'''}_1(x)=(3x-x^3)p_1(x).
\end{gather*}
Since $|p^{'''}_1(x)|$ is even and  $\int p^{'''}_1(x)dx=0$, we obtain
\begin{gather*}   %                                                                                     \label{g25}
\int |p_1^{'''}(x)|dx=2\int (p^{'''}_1(x))^{+}dx=4\int\limits^{\infty}_0(p^{'''}_1(x))^{+}dx
=4\int\limits^{\sqrt{3}}_0 (3x-x^3)p_1(x)dx
\\=4\int\limits^{\sqrt{3}}_0\!\! p_1'''(x)dx
=4p_1''(\sqrt{3})-4p_1''(0)
%=8p_1(\sqrt{3})+4p_1(0)
=\frac{4}{\sqrt{2\pi}}(1+\frac2{e^{3/2}})<2.4.
\end{gather*}
Thus, substituting this fact  into (\ref{g20}), we obtain
the desired inequality   (\ref{g16}).
\hfill\fbox

\section
{Proofs of Theorem \ref{thm2}  and  it's  corollaries }
In all cases when we are given $n$ values,
say $x_1,\dots,x_n$, denoted by the same letter with indices,
we agree  to use everywhere  vector notation:\ %s as follows
$\vec x_k=(x_1,\dots,x_k)$, for $k=1,\dots,n$.

\medskip\noindent{\bf 4.1.  Representations for broken lines.}
Define
$$\zeta_i:=\big(X_i-\mu_{i})\big/\big(\sqrt{n}\sigma_{i}\big),  \qquad  i=1,2,\ldots, n,$$
where the values $\mu_{i}$ and $\sigma_{i}$ are continuous functions of $\vec X_{i-1}=(X_1, X_2, \ldots, X_{i-1})$ which will be specified later.
Then   the random broken line $S_n$ with the values (\ref{a27}) at key points can be rewritten  as
\begin{gather}                                                                                \label{f2}
S_n(t)=\sum_{i=1}^{n}\zeta_{i}e_{i}(t), \  \   \  t\in [0,1],
\end{gather}
where growing functions $e_i$ have the following form:
\begin{gather}                                                                                               \label{f3}
e_i(t) = e_{i,n}(t):=
\begin{cases}
0,    \ \text{for}\ \ t\leq \frac{i-1}{n},\\
 nt-(i-1), \ \text{for}\ \  t\in [\frac{i-1}{n}, \frac{i}{n}] ,\\
1,    \ \text{for}\ \ t\geq \frac{i}{n}.
\end{cases}
\end{gather}
Here $i=0,1,\dots, n+1$. Note that
$e_{0,n}(t)=1$ and  $e_{n+1,n }(t)=0$   for all $t\in [0,1].$

Let $B=B(\cdot)$ be a standard  Wiener process, i.e.,
${\bf E}B(t)=0$ and ${\bf E}B^2(t)=t$ for all  $t\in  [0,1].$
Then for each $n=1,2,\dots$ and  all $t\in  [0,1]$ define
\begin{gather*}           %                                                                     \label{f7}
B_n(t):=\sum_{i=1}^{n}\theta_{i}e_{i}(t)
\quad\text{with}\quad
 \theta_i := \theta_{i,n}=B(\frac{i}{n})-B(\frac{i-1}{n}).
%%, \ \   \forall  t\in[0,1].
\end{gather*}
Thus,  $B_n=B_n(\cdot)$  is  also  a  random broken line.

%%\noindent{\bf 4.2. Preliminary smoothing.}\
\medskip\noindent{\bf 4.2.   Main idea.}
We are going to describe a method how it is possible
%%in sublinear expectation from (\ref{f12})
to change in (\ref{f2}) variables $\zeta_n, \zeta_{n-1}, \ldots, \zeta_2, \zeta_1$ one by one into $\theta_n, \theta_{n-1},$ $ \ldots, $  $\theta_2,$ $ \theta_1$, and to obtain (\ref{b11}) as a result.
Fix an integer $n\ge1$, real numbers $\rho,r>0$ and an arbitrary set $D\subset \mathbb{C}[0,1].$
Introduce a sequence $\nu,\nu_1,\nu_2,\dots$ of i.i.d. classical random variables with density
$g_1(\cdot)$ from (\ref{g1}), which we suppose to be independent from the Wiener process $B$.
Underline that $|\nu|<1$ and $|\nu_i|<1$ a.s. for all $i$.
We have in mind that almost surely
\begin{gather}                                                                                \label{f10}
|r\nu|<r\quad\text{and}\quad
\parallel\sum_{i=1}^n \frac{\rho}{n}\nu_i e_i\parallel
=\max_{t\in [0,1]}\Big| \sum_{i=1}^n \frac{\rho}{n}\nu_i e_i(t)\Big|
\le\sum_{i=1}^n \frac{\rho}{n}=\rho
\end{gather}
for some negligible $\rho>0$ which will be chosen later.
We have from (\ref{f10}) that
\begin{gather}                                                                                \label{f11}
h_{n+1}(\vec x_n):={\bf P}\big(\sum_{i=1}^n (x_i+\frac{\rho}{n}\nu_i)e_i\in D^{\rho}    \big)
\ge {\bf I}( \sum_{i=1}^n x_i e_i\in D).
%,\  \  \vec x_n=(x_1, x_2, \ldots, x_n).
\end{gather}

In the next subsection we construct a special sequence of sets $D_{k}$ which have, among others, the following property
\begin{gather}                                                                                                \label{f12}
D^{\rho}= D_{ n+1}\subset \ldots \subset D_{ k+1}\subset D_{k}\subset \ldots  \subset D_0\subset D^{\rho+r},\quad 1\le k\le n.
\end{gather}
In particular, it follows from (\ref{f10}) and (\ref{f12}) that
\begin{gather}                                                                                                \label{f13}
h_0:={\bf P}(B_n+r\nu e_0+\sum_{i=1}^n \frac{\rho}{n}\nu_ie_i\in  D_0)\le {\bf P}(B_n\in  D^{2\rho+2r}).
\end{gather}
Definition (\ref{a26}) and (\ref{f11}) imply that
$\overline{{\bf P}}(S_n\in D)\le  \ooo h_{n+1}(\vec\zeta_n)$.
Using also  (\ref{f13}) we obtain:
\begin{gather*}     %                                                                                           \label{f49}
\Delta:=\overline{{\bf P}}(S_n\in D)-{\bf P}(B_n\in D^{2\rho+2r})
%\le \overline{E}h_{n+1}(\vec\zeta_n)-h_0
\le \ooo h_{n+1}(\vec\zeta_n)-h_0.
\end{gather*}
Thus, for any continuous and bounded functions $h_1,\dots, h_n$
\begin{gather}                                                                                                 \label{f15}
\Delta\le\sum_{k=1}^{n+1} \Delta_k\   \    \text{ with } \  \
\Delta_k=\ooo h_k(\vec\zeta_k)-\ooo h_{k-1}(\vec\zeta_{k-1}).
\end{gather}

Values $\Delta_k$ will be estimated in subsection 4.6
with $h_k$ from (\ref{f17}) below, and
with special $\mu_{k}$  and  $\sigma_{k}$ which will be
constructed   in subsection 4.5.

\medskip\noindent{\bf 4.3.   Construction of sets $D_k$.}\
Similar to Sakhanenko  \cite{S88}  and  \cite{S00}, we  introduce a sequence of sets $D_{k}$ which we will interpret  as the  sequence of enlargements of the set $D$. We first define $D_{ n+1}=D^{\rho}$. After that for each  $k=n,n-1,\dots ,2,1,0$ we introduce sets $D_k$ by the following rule:
\begin{gather}      \label{f16}
D_{k}=D_k(r,\rho):=\big\{\tilde x=\tilde y+z_k(e_k-e_{k+1})\in \mathbb{C}[0,1] : \tilde y\in D_{ k+1}, |z_k|<r\big\}.
\end{gather}
Now for $k=1,\dots,n$ introduce  functions
\begin{gather}                                                                                                \label{f17}
h_{k}(\vec x_k)
%%=h_{k}(\vec x_{k-1},x_k)
:=
{\bf P}\Big(\sum_{i=1}^{k}(x_i+\frac{\rho}{n}\nu_i)e_i+r\nu e_k+\sum_{j=k+1}^n( \theta_j+\frac{\rho}{n}\nu_j) e_j\in D_{k}\Big),
\end{gather}
where we assume that $\sum_\emptyset=0$.
Since $e_n=e_n-e_{n+1}$, we have from definitions (\ref{f11}) and (\ref{f16}) of $D_{n+1}$ and  $D_n$ that
\begin{gather}                                                                                                \label{f19}
h_{n+1}(\vec x_n)\le
h_{n}(\vec x_n)
={\bf P}\big(r\nu e_n+\sum_{i=1}^n (x_i+\frac{\rho}{n}\nu_i)e_i\in D_n    \big);
\   \    \text{ so that } \  \ \Delta_{n+1}\le0.
\end{gather}
%%So that $\Delta_{n+1}\le0$.

\begin{lem}                                                                                                                \label{LL-6}
For each $r>0$ and $\rho>0$ inclusion of sets (\ref{f12}) takes place for sets $D_k$ defined in (\ref{f16}).
Moreover for $k=1,\dots,n$
\begin{gather}\label{f20}
{\bf E}h_k(\vec x_{k-1}, \theta_k)\le  h_{k-1}(\vec x_{k-1}).
\end{gather}
\end{lem}
{\bf    Proof. }
From   the  definition  (\ref{f16})   we  have
\begin{gather*}      %                                                                                        \label{f20}
D_0=\big\{\tilde   x=\tilde y+\tilde z: \tilde   y\in D^{\rho}, \tilde z=\sum_{k=0}^nz_k(e_k-e_{k+1}), \max|z_k|<r
\big\}.
\end{gather*}
In particular,
$$\parallel\tilde  z\parallel=\max_{t\in [0,1]}\big|\tilde z(t)\big|\le r\max_{t\in[0,1]}\sum_{k=0}^n\big|e_k(t)-e_{k+1}(t)\big|.$$
But $e_k(t)=e_0(t-k/n)$  for  all  $t$   as   it  follows   from  definition (\ref{f3})  of  non-decreasing functions~$e_k.$ Hence
$$\sum_{k=0}^n\big|e_k(t)-e_{k+1}(t)\big|=\sum_{k=0}^n\Big[e_0\big(t-\frac kn\big)-e_0\big(t-\frac{k+1}n\big)\Big]=e_0(t)-e_0\big(t-\frac{n+1}n\big)\le 1.$$
Thus  $\parallel\tilde  z \parallel<r$  and  we  have  from (\ref{a16})  that
$D_0\subset D^{\rho+r}.$

All other inclusions in (\ref{f12})   follow immediately from  definitions (\ref{f16}).

Now we are going to prove the main property (\ref{f20}) of sets $D_k$ introduced in  (\ref{f16}).
From   definition  (\ref{f17}) we have:
\begin{gather}\label{f20+}
{\bf E}h_k(\vec x_{k-1}, \theta_k)
={\bf P}\Big(\sum_{i=1}^{k-1}(x_i+\frac{\rho}{n}\nu_i)e_i+r\nu e_k+\sum_{j=k}^n( \theta_j+\frac{\rho}{n}\nu_j) e_j\in D_{k}\Big).
\end{gather}
But if \ $\tilde x+r\nu e_k\in D_k$ for some $\tilde x\in \mathbb{C}[0,1]$,
then  it follows from  the  definition  (\ref{f16}) of  $D_{k-1}$ with
$z_{k-1}=r\nu$ that
$$
\tilde x+r\nu e_{k-1}=\tilde x+r\nu e_k+z_{k-1}( e_{k-1}-e_k)\in D_{k-1}.
$$
This fact and (\ref{f20+}) together with (\ref{f17}) imply that
\begin{gather*}%\label{f20+}
{\bf E}h_k(\vec x_{k-1}, \theta_k)
={\bf P}(\tilde x+r\nu e_k\in D_{k})\le{\bf P}(\tilde x+r\nu e_{k-1}\in D_{k-1})
=h_{k-1}(\vec x_{k-1})
\end{gather*}
with evident random  $\tilde x$.
Thus,~(\ref{f20}) is proved.

\hfill\fbox

%\medskip
\noindent{\bf 4.4.   Using Taylor formula.}\
All broken lines under probability in (\ref{f17}) has summand
 $\sum_{i=1}^n \frac{\rho}{n}\nu_i e_i(t).$
We are going to show that it makes all functions $h_{k}$ to be sufficiently smooth.

\begin{lem}                                                                                                              \label{LL-5}
For each $k=1,\dots,n$  function $h_{k}(\vec x_k)=h_{k}(\vec x_{k-1},x_k)$
is twice continuously differentiable  with respect to its last variable $x_k$.
In addition, %%%in this case
the following functions
\begin{gather}                                                                                                \label{f22}
h_{k}'(\vec x_{k-1},0)=\frac{\partial}{ \partial x}h_{k}(\vec x_{k-1},x)\big|_{x=0}
\quad     \text{and }\ \
h_{k}''(\vec x_{k-1},0)=\frac{\partial^2}{\partial x^2}h_{k}(\vec x_{k-1},x)\big|_{x=0}
\end{gather}
are uniformly continuous functions of $\vec x_{k-1}$.
%%where the density $g_r(\cdot)$ was introduced in (\ref{c31}).
\end{lem}
{\bf  Proof. }
Consider  the  set
$D^{*}_k=\big\{\vec x_n=(\vec x_{k-1}, x_k, \vec x_{k,n}): \sum_{i=1}^n x_ie_i\in D_k    \big\}.$
It  is  easy  to  see  that
\begin{gather}
h_{k}(\vec x_k)=P\big( (\vec x_{k-1}+\frac{\rho}{n}\vec\nu_{k-1}, x_k+r\nu+\frac{\rho}{n}\nu_k, \vec\theta_{k,n}+\frac{\rho}{n}\vec\nu_{k,n})\in D^{*}_k   \big)\nonumber\\
\label{f23}
=\!\int\!\! dz\!\underbrace{\int\!...\! \int\!}_{D^{*}_k}\prod_{i=1}^{k-1} g(y_i-x_i)g(y_k)g_r(z-x_k-y_k)
\tilde{g}(\vec y_{k,n})
%%\prod_{j=k+1}^n \tilde{g}(y_j)
dy_1\ldots dy_n ,
\end{gather}
where  $g(\cdot)=g_{\rho/n}(\cdot)$  is  the  density  of  $\frac{\rho}{n}\nu_i$   for  all  $i,$  $g_r(\cdot)$  is  the  density  of  $r\nu$
and  $\tilde{g}(\cdot)$  is  the  joint density  of  $\vec\theta_{k,n}+\frac{\rho}{n}\vec\nu_{k,n}$.
%%$\theta_j+\frac{\rho}{n}\nu_j$   for  each  $j.$
Using (\ref{f23}), it  is  not  difficult to  verify  that for  all  $l<k$
$$\frac{\partial^3 h_k(\vec x_k)}{\partial x_l \partial x^2_k}
=\!\int\!\! dz\!\underbrace{\int\!...\! \int\!}_{D^{*}_k}
\prod_{
\substack{i=1 \\ i\neq l
}}^{k-1}g(y_i-x_i)g'\!(y_l-x_l)g(y_k)g''_r(z-x_k-y_k)
%%\prod_{j=k+1}^n
\tilde{g}(\vec y_{k,n})dy_1\!... dy_n.$$
Thus, with $z_i=y_i-x_i$ and $z_k=z-x_k-y_k$ we obtain that
\begin{gather} \label{f24}
\Big|  \frac{\partial^3 h_k(\vec x_k)}{\partial x_l \partial x^2_k} \Big|\le
\int \!\!\int\!\! \int \big| g^{'}(y_l-x_l) g^{''}_r(z-x_k-y_k) \big|g(y_k)dy_l dzdy_k\\  \nonumber
=\int |g^{''}_r(z_k)|dz_k \cdot \int |g^{'}(z_l)|dz_l
%%\nonumber\\&=&
=\frac{1}{r^2} \int |g^{''}_1(x)|dx \cdot \frac{n}{\rho}\int |g^{'}_1(x)|dx<\infty.
\end{gather}
%where  $z_i=y_i-x_i, \  \  z_k=z-x_k-y_k.$

So,  from  (\ref{f24})   we  have  that  $h^{''}_k(\vec x_{k-1}, 0)$   is  uniformly  continuous  because  all  its  partial  derivatives  are  uniformly bounded.
  The   similar   arguments  are true  also   for    function
 $h_k(\vec x_k)$ and for $h^{'}_k(\vec x_{k-1},0)$ from  (\ref{f22}).

\hfill\fbox

Introduce  an error  term  of the Taylor  expansion  of  function $h_{k}(\vec x_{k-1},x_k)$
with respect to its last variable $x_k=x$ at $x_k=x=0$:
\begin{gather}                                                                                      \label{f28}
\delta_{k}(y|\vec x_{k-1})=\delta(0, y|h) \quad\text{for}\quad
h(\cdot)=h_{k}(\vec x_{k-1},\cdot),\quad k=1,\dots,n,
\end{gather}
where the function $\delta(x,y|h)$ was defined in (\ref{c7}).

\begin{lem}                                                                                                              \label{LL-7}
If  $ p\in[2,3],$ then for each $k=1,\dots,n$ and all $y\in \mathbb{R}$
\begin{gather}                                                                                                \label{f26}
\delta_{k}( y):=\sup_{\vec x_{k-1}}
\big|\delta_{k}( y|\vec x_{k-1})\big|\le \overline{\varepsilon}_{p,r}(y)
:=16 \min\Big\{\frac{ |y|^p}{(2r)^p}, \frac{ y^2}{(2r)^2}\Big\}.
%,\quad p\in[2,3].
\end{gather}

\end{lem}
{\bf  Proof. }
Introduce  set:
$$A_k(\vec y_{k-1}, y, \vec y_{k,n}))
=\big\{x: \sum_{i=1}^{k-1}y_ie_i+(x+y)e_k+\sum_{j=k+1}^ny_je_j\in D^k    \big\}.$$
Then from  (\ref{f17})  we  have
$$
h(x)
%=h_k(\vec x_k)
=h_k(\vec x_{k-1}, x)
={\bf P}\big(x+r\nu\in A_k(\vec x_{k-1}+\frac{\rho}{n}\vec\nu_{k-1},\frac{\rho}{n}\nu_k, \vec\theta_{k,n}+\frac{\rho}{n}\vec\nu_{k,n})\big).
$$
Applying  Lemma \ref{LL-2} to the random set
$A_k\big(\vec x_{k-1}+\frac{\rho}{n}\vec\nu_{k-1},
 \frac{\rho}{n}\nu_k, \vec\theta_{k,n}+\frac{\rho}{n}\vec\nu_{k,n})\big),$
we obtain from (\ref{g6}) that
$\delta_{k}( y)\le \overline{\varepsilon}_{r}(y)$ with a function
$\overline{\varepsilon}_{r}(y)\le \overline{\varepsilon}_{p,r}(y)$ by (\ref{g2}).
\hfill\fbox

\noindent
{\bf 4.5.  Constructions   of $\mu_{k}$  and  $\sigma_{k}$.}\
Choose numbers $\varkappa >0$ and $\kappa >0$ so that
%%and $\rho$ in the following way:
\begin{gather}                                                                                      \label{f30}
\frac{\varkappa (\overline{\mu}-\underline{\mu})}{2\sqrt{n}\underline{\sigma}}
%%=\frac{\rho_0}{(\varepsilon\sqrt n)^3},\quad
+ \frac{\kappa }{2n}(\frac{\overline{\sigma}^2}{\underline{\sigma}^2}-1)
\le\rho_0\min_{k\le n}{{\bf E}}\Big[\delta_k\Big( \frac {Z}{\sqrt n}\Big)\Big]
%{(\varepsilon\sqrt n)^3}
,\quad
\rho_0=0.001,
\end{gather}
where $Z\sim N(0,1)$. After that define functions
\begin{gather}                                                                                      \label{f31}
\alpha_{k}(\vec x_{k-1}):=\alpha_{\varkappa }(0|h) \quad\text{for}\quad
h(\cdot)=h_{k}(\vec x_{k-1},\cdot),\quad k=2,\dots,n,
\end{gather}
where function $\alpha_{\varkappa }(0|h)$ was defined in (\ref{c1}) and (\ref{c2});
and we put
\begin{gather}                                                                                      \label{f32}
\beta_{k}(\vec x_{k-1})=\beta_{\varkappa ,\kappa }(0|h) \quad\text{for}\quad
h(\cdot)=h_{k}(\vec x_{k-1},\cdot),\quad k=2,\dots,n
\end{gather}
with function $\beta_{\varkappa ,\kappa }(0|h)$ introduced in (\ref{c3}) and (\ref{c4}).

\begin{lem}                                                                \label{LL-8}
For each $k=2,\dots,n$ functions $\alpha_{k}(\vec x_{k-1})$ and $\beta_{k}(\vec x_{k-1})$
%%%, defined by formulas (\ref{f31}) and (\ref{f32}),
are continuous functions of $\vec x_{k-1}\in \mathbb{R}^{k-1}$ and satisfy properties (\ref{c5+}). Moreover,
for all   $\vec x_{k-1}\in \mathbb{R}^{k-1}$
% and  $k=1,\dots,n$
\begin{gather}                                                                                      \label{f33}
H_k(\vec x_{k-1}):=
\ooo h_{k}(\vec x_{k-1}, \xi_k(\vec x_{k-1})/\sqrt n)\le h_{k-1}(\vec x_{k-1}) +\tau_{k,n},
\end{gather}
where $\xi_k(\vec x_{k-1}):=(X_0-\alpha_{k}(\vec x_{k-1}))/\beta_{k}(\vec x_{k-1})$ and
\begin{gather}                                                                                      \label{f34}
\tau_{k,n}= \sup_{\mu\in[\underline{\mu},\overline{\mu}]}\ooo\Big[\delta_k\Big(
\xi_0(\mu)/\sqrt n\Big)\Big]
+(1+\rho_0){{\bf E}}\Big[\delta_k\Big(Z/{\sqrt n}\Big)\Big].
\end{gather}
\end{lem}
{\bf  Proof. }
We  are  going  to  apply   Lemma \ref{LL-1}   at  $x=0$  to  auxiliary   function  $h(\cdot)$  used  in  (\ref{f28}), (\ref{f31})  and  (\ref{f32}).  Note  that in  this  case $\xi(0)=\xi_k(\vec x_{k-1})$  and  $Z=\sqrt{n}\theta_k$   has  standard normal  distribution.  Hence, using (\ref{f26}) and  (\ref{f30})  we  obtain  from  (\ref{c6})  with  $x=0$  that
\begin{gather}                                                                                     \label{f36}
\ooo h_k(\vec x_{k-1}, \frac{\xi_k(\vec x_{k-1})}{\sqrt{n}})\le {\bf E}h_k(\vec x_{k-1}, \theta_k)+\ooo \delta_k(\frac{\xi_k(\vec x_{k-1})}{\sqrt{n}})
+(1+\rho_0){\bf E}\delta_k(\theta_k).
\end{gather}
Note that $|\xi_k(\vec x_{k-1})|\le|\xi_0(\mu)|$ with
$\mu=\alpha_{k}(\vec x_{k-1})$ as it follows from (\ref{c5+}).
This fact together with (\ref{f36}) yields
\begin{gather}                                                                                      \label{f37}
H_k(\vec x_{k-1})=\ooo h_{k}(\vec x_{k-1}, \xi_k(\vec x_{k-1})/\sqrt n)\le
{\bf E}h_k(\vec x_{k-1}, \theta_k) +\tau_{k,n},
\end{gather}
Here we used also definitions from (\ref{f33}) and (\ref{f34}).

Therefore  (\ref{f33})  follows  from  (\ref{f37})  and  (\ref{f20}).
\hfill\fbox

Now  for  $k=1,2, \ldots, n$  we  are  going  to   define  functions  $\mu_k(\vec x_{k-1})$  and  $\sigma_k(\vec x_{k-1}).$  We  first  introduce
\begin{gather}\label{f41}
\mu_1=\mu_1(\vec x_0)=\underline{\mu} \  \   \text{and}\  \   \sigma_1=\sigma_1(\vec x_0)=\underline{\sigma}(\underline{\mu})>0.
\end{gather}
Note  that  if  for  some  $k>1$   all $\mu_i$  and  $\sigma_i$
%% functions  $\mu_i(\vec x_{i-1})$  and  $\sigma_i(\vec x_{i-1})$
 for  $i<k$  are  defined,  then  we  may  put
$$z_i(\vec x_i)={\big(x_i-\mu_i(\vec x_{i-1})\big)}\big/\sqrt{n}{\sigma_i(\vec x_{i-1})}, \   \   i=1,2,\ldots,k,$$
and  define  vector $\vec z_{k-1}(\vec x_{k-1})
%$\vec z_k(\vec x_k)
=(z_1(\vec x_1), \ldots, z_{k-1}(\vec x_{k-1})).$
At  last, after  that  we  may   introduce
\begin{gather}\label{f44}
\mu_k(\vec x_{k-1})=\alpha_k(\vec z_{k-1}(\vec x_{k-1})), \   \
\sigma_k(\vec x_{k-1})=\beta_k(\vec z_{k-1}(\vec x_{k-1})),\ \ k=2,\dots,n.
\end{gather}

 So, later on we use in our proof only
%\begin{gather*}                                      %               \label{f45}
$\zeta_k=z_k(\vec X_k)=(X_k-\mu_k)/\sqrt{n}\sigma_k$% \ \ \  \text{with} \   \
with $\mu_k=\mu_k(\vec X_{k-1})$ %\  \   \text{and} \   \
and $\sigma_k=\sigma_k(\vec X_{k-1})$
%%,  k=1,2,\ldots,n. \end{gather*}
defined in (\ref{f41}) and (\ref{f44}).

%%%\noindent{\bf 4.5.  Using independence under sublinear expectation.}\

\begin{lem}                                                                \label{LL-9}
If all values $\mu_k$ and $\sigma_k$
%, $k=1,\dots,n$,
are defined in (\ref{f41}) and (\ref{f44}), then
\begin{gather}                                                                                                \label{f47}
\Delta:=\overline{{\bf P}}(S_n\in D)-{\bf P}(B_n\in D^{2\rho+2r})\le \sum_{k=1}^n\tau_{k,n}.
%%p(2r\sqrt n).
%%\le\frac{5.2\gamma_3(2r\sqrt n)}{\sqrt{n}r^3}.
\end{gather}
In addition, for all $k=1,\dots,n$
\begin{gather}\label{f48}
\tau_{k,n}\le 16\Big\{\frac{\gamma_p(2r\sqrt{n})}{(2r\sqrt{n})^p}+(1+\rho_0)
{\bf E} \min\Big\{ \frac {|Z|^p}{(2r\sqrt n)^p},\frac {Z^2}{(2r\sqrt n)^2}\Big\}
 \Big\}.
\end{gather}
\end{lem}

{\bf  Proof. }
Introduce  function:
\begin{gather}\label{f51}
\overline{h}_k(\vec x_{k-1}, x_k)=\overline{h}_k(\vec x_k)=h_k(\vec z_k(\vec x_k))=h_k(\vec z_{k-1}(\vec x_{k-1}),
\frac{x_k-\mu_k(\vec x_{k-1})}{\sqrt{n}\sigma_k(\vec x_{k-1})}).
\end{gather}
Remind  that  function  $\overline{h}_k(\cdot)$  is  continuous   and  bounded,  whereas  $X_k$   and  $X_0$   are  identically  distributed.  Hence
\begin{gather}\label{f52}
\ooo\overline{h}_k(\vec x_{k-1}, X_k)=\ooo h_k(\vec z_{k-1}(\vec x_{k-1}), \frac{X_0-\mu_k(\vec x_{k-1})}{\sqrt{n}\sigma_k(\vec x_{k-1})})=H_k(\vec z_{k-1}(\vec x_{k-1})).
\end{gather}
In addition   we   used here  definition  (\ref{f33})  of  $H_k(\cdot).$
But  $X_k$   is  independent  to  $\vec X_{k-1}.$  Hence  using  definition  (\ref{a20}) of  independence  with  $\psi(\cdot)=\overline{h}_k(\cdot)$  we  obtain:
\begin{gather}\label{f53}
\ooo   \overline{h}_k(\vec X_k)=\ooo\big[ \ooo[h_k(\vec x_{k-1}, X_k)]
  \big|_{\vec x_{k-1}=\vec X_{k-1}}\big]
  =\ooo H_k(\vec z_{k-1}(\vec X_{k-1}))
  =\ooo H_k(\vec\zeta_{k-1}).
\end{gather}

From  (\ref{f51})--(\ref{f53})  we  have :
$$\ooo h_k(\vec\zeta_k)=\ooo\overline{h}_k(\vec X_k)=\ooo H_k(\vec\zeta_{k-1})
\le \ooo(h_{k-1}(\vec\zeta_{k-1}))+\tau_{k,n}
= \ooo h_{k-1}(\vec\zeta_{k-1})+\tau_{k,n}
.$$
Here  we  also  used  inequality (\ref{f33})  with $\vec x_{k-1}=\vec\zeta_{k-1}.$
So  we  obtained  that  $\Delta_k\le  \tau_{k,n}.$
Substituting  the  last  estimate  into  (\ref{f15}), we obtain (\ref{f47})
if only we remind that $\Delta_{n+1}\le 0$ by (\ref{f19}).

Inequality (\ref{f48}) follows immediately from  (\ref{f26}) and (\ref{f34}).
\hfill\fbox

\begin{rem}                                                        \label{R41}
If $X_1,X_2,\dots$ are classical random variables then all arguments in the present proof of Theorem \ref{thm2} remains valid with $\rho=\varkappa =\kappa =0$ for all measurable sets $D$. As a result from (\ref{f47})
with $r=\varepsilon/2$ we obtain that
\begin{gather}                                                                                                \label{f72}
\overline{{\bf P}}(S_n\in D)-{\bf P}(B_n\in D^{\varepsilon})\le
\frac {16}{\varepsilon^2}{\bf E} \min\Big\{ \frac {|\xi_1|^3}{\varepsilon\sqrt n},\xi_1^2\Big\}
+\frac {16}{\varepsilon^2}{\bf E} \min\Big\{ \frac {|Z|^3}{\varepsilon\sqrt n},Z^2\Big\},
\end{gather}
where $Z$ has standard normal distribution.
This estimate was first obtained in Corollary 3.1 in Sakhanenko  \cite{S00}, using similar, but simplified, arguments.

Estimate (\ref{f72}) may be treated as the first estimate in FCLT with explicit constants. But
(\ref{f72}) is not a "real"  estimate in FCLT because it is not an estimate for the value:
\begin{gather} \label{f73}
\Delta_+(\varepsilon):=\overline{{\bf P}}(S_n \in D)-{\bf P}(B\in D^{\varepsilon})
\end{gather}
which we are going to estimate in the present paper.
\end{rem}

\noindent{\bf 4.6. Proof of Theorem {\ref{thm2}}. }\
We need two auxiliary lemmas.
\begin{lem}                                                                \label{LL-11}
If $p\ge2$ then $\gamma_p(C)\ge K_p(C):=\min\big\{1, 2C^{p-2}/p\big\}$  for any $C>0$.
\end{lem}
{\bf  Proof. }  Consider  function  $f_0(x)=\int\nolimits^{|x|}_0 (p/2) \min\big\{y^{p/2-1}, 1 \big\}dy.$  It is  even and  convex  because  it's derivative  $f^{'}_0(x)$   is  non-decreasing  for  $x\ge 0.$  By  Jensen inequality  for  sublinear  expectations   we  have  for any  $\mu\in [\underline{\mu}, \overline{\mu}]$  that
\begin{gather}\label{f74}
\ooo f_0(\xi^2_0)\ge f_0(\ooo [\xi^2_0])=f_0(\overline{\sigma}^2(\mu)/\underline{\sigma}^2(\mu))\ge f_0(1)=1,
\end{gather}
where $\xi_0=(X_0-\mu)/\underline{\sigma}(\mu)$
and $f_0(1)=\int\nolimits^{1}_0 (p/2)y^{p/2-1}dy=1.$

Next  consider   function
$f(x)=\min\big\{|x|^{p/2}, C^{p-2}|x|    \big\}.$
It  is  easy  to  see  that
%%for  all  $x\ge 0$
\begin{gather*}
\forall\ x\ge0\quad f^{'}(x)=
\begin{cases}
px^{p/2-1}/2,    \ \text{if}\ \ x< C^2\\
C^{p-2}, \ \text{if}\ \  x> C^2
\end{cases}
\ge K_p(C)f^{'}_0(x).
\end{gather*}
Hence, $f(x)=f(|x|)\ge K_p(C)f_0(|x|)=K_p(C)f_0(x).$  From  this  fact, (\ref{f74}), (\ref{b2})  and  (\ref{b3})  we obtain  that
$$\gamma_p(C)\ge \gamma_p(u, C)=\ooo f(\xi^2_0)\ge K_p(C)\ooo f_0(\xi^2_0)\ge K_p(C).$$
\hfill\fbox

\begin{lem}                                                                        \label{LL-10}
For any $b>0$
\begin{gather}                                                                  \label{f75}
\Delta_n^*:={\bf P}(B_n \in D^{2\rho+2r})-{\bf P}(B\in D^{2\rho+2r+b})\le {\bf P}(\|B_n-B\|\geq b)
\leq \frac{n{\bf E}|Z|^{p}}{(\sqrt nb)^{p}}.
\end{gather}

\end{lem}
{\bf Proof. }
It is well known that function $|x|^p$ is convex for $p\ge1$ and $B(t)$ is a martingale.
So, by maximal   inequality we have
\begin{gather*}
P_n^*(b):={\bf P}(\sup_{t\in [0,{1}/{n}]}|B(t)|\geq  b)\le\frac{{\bf E}|B(1/n)|^{p}}{b^{p}}=\frac{{\bf E}|B(1)|^{p}}{(\sqrt nb)^{p}}=
\frac{{\bf E}|Z|^{p}}{(\sqrt nb)^{p}}.
\end{gather*}
Note that $B(k/n)=B_n(k/n)$ for $k=0,1,\dots,n$. Hence
\begin{eqnarray*}
{\bf P}(\|B_n-B\|\geq b)&=& {\bf P}(\max _{t\in [0,1]}|B_n(t)-B(t)|\ge b)\\
&\leq&\sum_{k=1}^{n}{\bf P}(\sup_{t\in [\frac{k-1}{n},\frac{k}{n}]}|B_n(t)-B(\frac{k-1}{n})|\geq  b)=nP_n^*(b).
\end{eqnarray*}
Since the first inequality in (\ref{f75}) is evident,  the  result  is  proved.

\hfill\fbox

Thus,  from (\ref{f47}), (\ref{f75}), (\ref{f73})   and  (\ref{f48}) with $b=r\le\varepsilon/2$ we have:
\begin{gather} \label{f77}
%%\overline{P}(S_n \in D)-P(B\in D^{3r+2\rho})
\Delta_+(3r+2\rho)=\Delta+\Delta_n^*
\le
\frac{n}{(r\sqrt n)^p}\Big[\frac{16\gamma_p(\varepsilon\sqrt n)}{2^p}+
\frac{16(1+\rho_0){\bf E}|Z|^{p}}{2^p}+{\bf E}|Z|^{p}\Big].
\end{gather}
Choose $\rho=\rho_0\varepsilon/2$ and $r=(\varepsilon-\rho_0\varepsilon)/3$ in (\ref{f77}) and assume that
$(\varepsilon\sqrt n)^{p-2}\ge p/2$. As a result we obtain, by Lemma \ref{LL-11}, the desirable inequality (\ref{b11}), but with $C_5(p)$ instead of $C_2(p)$, where
\begin{gather*}    %                                                                                            \label{f78}
C_5(p):=
\frac{3^p}{(1-\rho_0)^p}\Big[\frac{16}{2^p}+
\frac{16(1+\rho_0){\bf E}|Z|^{p}}{2^p}+{\bf E}|Z|^{p}\Big].
\end{gather*}
Since ${\bf E}|Z|^{p}\le({\bf E}|Z|^3)^{p/3}=({4}/{\sqrt{2\pi}})^{p/3}<(1.6)^{p/3}$,
it is not difficult to verify that $C_5(p)\le C_5(3)<184 $ and that $C_5(p)/(4.7)^p\le C_5(2)/(4.7)^2<4.7 .$
So, $C_5(p)\le C_2(p)$ for all $p\in[2,3]$.

Thus,  inequality (\ref{b11})  is  proved when   $C^{p-2}\ge p/2$  with
$C:=\varepsilon\sqrt{n}$. Consider now the case when $0<C^{p-2}=(\varepsilon\sqrt{n})^{p-2}< p/2$. Then $C^2\le ep/2$ and
%(\ref{b11})  also takes place because
$\gamma_p(C)\ge 2C^{p-2}/p$
by Lemma \ref{LL-11}.  Hence, in this case
the right hand side in (\ref{b11})  has the following property
\begin{gather*}       %                                                                                                    \label{f79}
\frac{C_2(p)\gamma_p(\varepsilon\sqrt{n})}{n^{(p-2)/2}\varepsilon^{p}}
=\frac{nC_2(p)\gamma_p(C)}{C^{p}}\ge\frac{2nC_2(p)C^{p-2}}{pC^{p}}=\frac{2nC_2(p)}{pC^{2}}\ge \frac{C_2(2)}e>1.
\end{gather*}
But the left hand side in (\ref{b11}) does not exceed 1 in view of (\ref{a26}).

So, Theorem \ref{thm2} is proved in all cases.

\hfill\fbox

\noindent{\bf 4.7. Proofs of Corollaries {\ref{cor1+}} and {\ref{cor2+}}. }
We are going to use the assertion of Theorem \ref{thm2} with
$\varepsilon=\min\big\{\tilde{\Pi}_{n,p},1\big\}.$ In this case
$\gamma_p(\varepsilon\sqrt{n})\le  \gamma_p(\sqrt{n}).$
Hence, for the chosen $\varepsilon$ we have from (\ref{b11}) that
$\tilde \Pi_n(\varepsilon,p)\le\varepsilon$.
Thus, by definition (\ref{b15}), the chosen $\varepsilon$
is an estimate for $\overline{\Pi}(S_n, B)$, and (\ref{b17}) follows.

Now we are going to apply Theorem \ref{thm2} with the same
$\varepsilon$
%%\noindent{\bf 4.8. Proof of Corollary {\ref{cor2+}}. }
to set  $D=D_a=\big\{x: f(x)\le a      \big\}.$  We obtain
\begin{gather}                                                                                                           \label{f81}
\overline{{\bf P}}(f(S_n)\le a)=\overline{{\bf P}}(S_n \in D_a)
\le  {\bf P}(B \in D^{\varepsilon}_a)+\tilde \Pi_n(\varepsilon,p)
\le {\bf P}(B \in D^{\varepsilon}_a)+\varepsilon.
\end{gather}
Since  $f(\cdot)$  satisfies   condition   (\ref{a17}), we  have that
$D^{\varepsilon}_a\subset D_{a+L\varepsilon}.$  Thus,  by
 condition (\ref{a18}) we   have
\begin{gather}                                                                                                           \label{f82}
 {\bf P}(B\in D^{\varepsilon}_a)
\le {\bf P}(B\in D_{a+L\varepsilon})= {\bf P}(f(B)\le a+L\varepsilon)\le  {\bf P}(f(B)\le a)+KL\varepsilon.
\end{gather}
Now (\ref{b19}) follows from (\ref{f81}) and (\ref{f82}).

\section{Proof of Theorem \ref{thm1CLT} } %% and  it's  corollaries }

%\textcolor{blue}{Continue, please, all that I have not done till the end}

\noindent{\bf 5.1. Preliminary considerations.}\
 We fix an integer $n\ge1$, a real number $\varepsilon>0$ and a  set $A\subset \mathbb{R}$. Introduce into consideration classical random variable $\nu$ with density
$g_1(\cdot)$ from (\ref{g1}), which we suppose to be independent from the Wiener process $B$.
Underline that $|\nu|<1$ a.s. and define
\begin{gather}\label{z1}
\tilde h_{k}(x):={\bf P}(x+\varepsilon\nu/2+B(1)-B(k/n)\in A^{\varepsilon/2}),\quad k=0,1,\dots,n.
\end{gather}
 For  any  $x, y \in \mathbb{R}$  introduce an  error term of  the Taylor  expansion of  $\tilde{h}_{k}(x+y):$
\begin{gather}                                                                                      \label{z2}
\tilde \delta_k(x,y):=\tilde h_k(x+y)-\tilde h_k(x)-y\tilde  h_k'(x)-y^2\tilde h_k''(x)/2.
\end{gather}

\begin{lem}                                                                \label{LL-12}
If $1\le k<n$ then for all $y\in \mathbb{R}$
\begin{gather}                                                                                                \label{z5}
\tilde \delta_{k}( y):=\sup_x\big|\tilde \delta_{k}(x, y)\big|\le0.4|y|^3(1-k/n)^{-3/2}.
%%\frac{2|y|^3\sqrt{n}^3}{5\sqrt{n-k}^3}
\end{gather}
\end{lem}
{\bf Proof.}
Note that random variable $Z=(B(1)-B(k/n))/b$ is independent of $\nu$ and has standard normal distribution when $b=\sqrt{1-k/n}.$
Hence, from definition (\ref{z1}) we have that
\begin{gather*}
\tilde h_{k}(x)={\bf P}\big(x+bZ \in A(\varepsilon\nu/2)\big)\quad     \text{ with } \  \
A(y)=A^{\varepsilon/2}-y.
%%k=0,1,\dots,n.
\end{gather*}
Thus, inequality (\ref{z5}) follows from the assertion (\ref{g16}) of Lemma \ref{LL-4}.

\hfill\fbox

\noindent{\bf 5.2.  Constructions   of $\mu_{k}$  and  $\sigma_{k}$.}\
Choose numbers $\varkappa >0$ and $\kappa >0$ so that
condition (\ref{f30}) be fulfilled with $\tilde \delta_k(\cdot)$ instead of $ \delta_k(\cdot)$.
After that define functions
\begin{gather}                                                                                      \label{z11}
\mu_{k}(x):=\alpha_{\varkappa }(x|h) \quad\text{for}\quad
h(\cdot)=\tilde h_{k}(\cdot),\quad k=2,\dots,n,
\end{gather}
where function $\alpha_{\varkappa }(x|h)$ was defined in (\ref{c1}) and (\ref{c2});
and we put
\begin{gather}                                                                                      \label{z12}
\sigma_{k}(x)=\beta_{\varkappa ,\kappa }(x|h) \quad\text{for}\quad
h(\cdot)=\tilde h_{k}(\cdot),\quad k=2,\dots,n
\end{gather}
with function $\beta_{\varkappa ,\kappa }(x|h)$ introduced in (\ref{c3}) and (\ref{c4}).
And  for  $k=1$  we  put
\begin{gather}                                                                              \label{z13}
\mu_1=\mu_1(x)=\underline{\mu} \  \   \text{and}\  \   \sigma_1=\sigma_1(x)=\underline{\sigma}(\underline{\mu})>0.
\end{gather}

 So, later on in our proof  we use only
$\zeta_k=(X_k-\mu_k)/\sigma_k$ with $\mu_k=\mu_k(\zeta_{k-1})$ and $\sigma_k=\sigma_k(\zeta_{k-1})$
defined in (\ref{z11}) -- (\ref{z13}).

\begin{lem}                                                                \label{LL-13}
If all values $\mu_k$ and $\sigma_k$
%, $k=1,\dots,n$,
are defined in (\ref{z11}) -- (\ref{z13}), then
\begin{gather}                                                                                                \label{z17}
\tilde\Delta:=\overline{{\bf P}}(W_n\in A)-{\bf P}(Z\in A^{\varepsilon})\le \sum_{k=1}^n\tilde \tau_{k,n},
\end{gather}
where
\begin{gather}                                                                                \label{z18}
\tilde\tau_{k,n}\le \frac {0.4C_7\gamma_3}{(n-k)^{3/2}}
 \ \ \  \text{with} \   \ \quad C_7:=1+(1+\rho_0){\bf E} |Z|^3<2.616,
\end{gather}
when $1\le k<n$. In addition,
\begin{gather}                                                                                \label{z19}
\tilde\tau_{k,n}\le \frac {16}{((\varepsilon-2\rho)\sqrt n)^p}\Big\{\gamma_p(\varepsilon\sqrt{n})+(1+\rho_0)
{\bf E} |Z|^p\Big\}\le \frac {16C_7\gamma_p}{((\varepsilon-2\rho)\sqrt n)^p},
\end{gather}
for all $k=1,\dots,n$, $p\in[2,3]$  and all $\rho\in(0,\varepsilon/2)$.
\end{lem}
{\bf  Proof. }
%%Introduce into consideration the following set
Consider a partial case of Theorem \ref{thm2} when
\begin{gather}                                                                                        \label{z21}
D=\big\{\tilde x(\cdot)\in \mathbb{C}[0,1]: \tilde x(1)\in A\big\}
 \ \ \  \text{and} \   \  r+\rho=\varepsilon/2.
\end{gather}
In this case $D^r=\big\{\tilde x(\cdot)\in \mathbb{C}[0,1]: \tilde x(1)\in A^r\big\}$ for all $r>0$ and hence
\begin{gather}                                                                                                \label{z22}
\overline{{\bf P}}(S_n\in D)=\overline{{\bf P}}(W_n\in A),
 \ \ \ \ \ %% \text{and} \   \
{\bf P}(B_n\in D^{\varepsilon})={\bf P}(Z\in A^{\varepsilon}),
\end{gather}
and $\Delta=\tilde \Delta$ as a result.
Moreover, it is not difficult to see from definitions (\ref{f16}), (\ref{f17}) and  (\ref{z1}) that for all $k=1,\dots,n$
\begin{gather*}      %%                                                                                          \label{z23}
D_k=D^{r+\rho}= D^{\varepsilon/2}
 \ \ \  \text{and} \   \
  h_k(\vec x_k)=\tilde h_k(x_1+\dots+x_k).
\end{gather*}
Using now definition (\ref{f28} ) with (\ref{z2}) and (\ref{f26} ) with (\ref{z5}) we obtain that
\begin{gather}                                                                                                \label{z24}
 \delta_k(\vec x_{k-1},y)=\tilde\delta_k(x_1+\dots+x_{k-1},y)
 \ \ \  \text{and} \   \
\delta_{k}( y)=\tilde \delta_{k}( y).
\end{gather}

Thus, under assumptions (\ref{z21}), from (\ref{z22}) and (\ref{z24}) we have that assertion (\ref{f47}) of Lemma \ref{LL-9} may be written in the form (\ref{z17}) with
\begin{gather}                                                                                      \label{z26}
\tilde\tau_{k,n}:= \sup_{\mu\in[\underline{\mu},\overline{\mu}]}\overline{\bf E}\Big[\tilde\delta_k\Big(\xi_0(\mu)/\sqrt n\Big)\Big]
+(1+\rho_0){{\bf E}}\Big[\tilde\delta_k\Big(Z/{\sqrt n}\Big)\Big]=\tau_{k,n}.
\end{gather}
Substituting (\ref{z5}) into (\ref{z26}) we obtain (\ref{z18}). At last, (\ref{z19}) is true as a corollary of inequality (\ref{f48}) obtained under more general assumptions.

\hfill\fbox

Note that estimates (\ref{z5}) and (\ref{z18}) are the only facts, used in the proof of Theorem {\ref{thm1CLT}},
which does not hold under more general conditions of Theorem {\ref{thm2}.

\medskip\noindent
{\bf 5.3.  Proof of Theorem {\ref{thm1CLT}}.}
If $(\varepsilon\sqrt n)^{p-2}\ge p/2$ then $\gamma_p(\varepsilon\sqrt n)\ge1$ by Lemma~\ref{LL-11}. In this case
substituting estimate (\ref{z19}) into (\ref{z17}) we find
inequality (\ref{b21}) with $C_4= {16C_7}/{(1-2\rho/\varepsilon)^3}\le42$
for $\rho$ sufficiently small.
When $(\varepsilon\sqrt n)^{p-2}<p/2$ we may repeat the arguments used at the end of the proof of Theorem \ref{thm2}
and obtain that (\ref{b21}) again holds because its right hand side is greater than 1.

Let $m$ be the minimal integer not smaller than $x_0=n\varepsilon^2/16$.
From (\ref{z18}) we  have:
\begin{gather}                                                                                                \label{z31}
 \sum_{k<n-m}\tilde \tau_{k,n}\le  \sum_{l>m} \frac {0.4C_7\gamma_3}{l^{3/2}}
< \int_{x>m}\frac {1.05\gamma_3}{x^{3/2}}dx<\frac {2.1\gamma_3}{x_0^{1/2}}=\frac {8.4\gamma_3}{\varepsilon\sqrt n}.
\end{gather}
On the other hand, from (\ref{z19}) with $p=3$ for sufficiently small $\rho>0$ we obtain:
\begin{gather}                                                                                                \label{z32}
 \sum_{k\ge n-m}^n\tilde \tau_{k,n}\le \frac {(m+1)16C_7\gamma_3}{((\varepsilon-2\rho)\sqrt n)^3}
\le (x_0+2)\frac {16\cdot2.7\gamma_3}{(\varepsilon\sqrt n)^3}=\frac {2.7\gamma_3}{\varepsilon\sqrt n}
+2\frac {16\cdot2.7\gamma_3}{(\varepsilon\sqrt n)^3}.
\end{gather}
Substituting estimates (\ref{z31}) and (\ref{z32}) into (\ref{z17}) we obtain (\ref{b22}) in the case when $\varepsilon\sqrt n\ge12$.
But if $\varepsilon\sqrt n<12$, then right hand side in (\ref{b22}) greater than 1. So, (\ref{b22}) takes place in all cases.

%\hfill\fbox

Thus, all assertions of the paper are proved.

\bigskip

{ \noindent {\bf\large  Acknowledgments} \quad  The   authors  are  grateful  to  the   anonymous referee
for his/her thorough review and highly appreciate his/her comments,
which inspired us  to improve the results and their presentation.

}

\end{document}